\documentclass[11pt]{article}

\usepackage{bridges}

\usepackage{pinlabel}
\usepackage{caption}
\usepackage{subcaption}
\usepackage{color}
\usepackage{amsmath}
\usepackage{amssymb}
\usepackage{hyperref}

\newcommand{\RN}[1]{%
  \textup{\uppercase\expandafter{\romannumeral#1}}%
}
\newtheorem{theorem}{\bf Theorem}[section]

\newtheorem{lemma}{\bf Lemma}[section]

\newtheorem{definition}{\bf Definition}[section]

\begin{document}

\title{A chiral family of triply-periodic minimal surfaces from the quartz network}

\author{
Shashank G. Markande $^{1}$, Matthias Saba$^{2}$, G. E. Schr\"{o}der-Turk$^{3}$ and Elisabetta A. Matsumoto $^{1*}$\\
\small $^{1}$School of Physics, Georgia Tech, Atlanta, GA 30318, USA\\
\small$^{2}$ Department of Physics, Imperial College, London, UK\\
\small$^{3}$School of Engineering and IT, Mathematics and Statistics, Murdoch University, Murdoch, Australia\\
\small$^{*}$ Author for correspondence: \href{mailto:sabetta@gatech.edu}{sabetta@gatech.edu}}





\maketitle

\begin{abstract}
We describe a new family of triply-periodic minimal surfaces with hexagonal symmetry, related to the quartz ({qtz}) and its dual (the {qzd} net). We provide a solution to the period problem and provide a parametrisation of these surfaces, that are not in the regular class, by the Weierstrass-Enneper formalism. We identified this analytical description of the surface by generating an area-minimising mesh interface from a pair of dual graphs (qtz $\&$ qzd) using the generalised Voronoi construction of De Campo, Hyde and colleagues, followed by numerical identification of the flat point structure. This mechanism is not restricted to the specific pair of dual graphs, and should be applicable to a broader set of possible dual graph topologies and their corresponding minimal surfaces, if existent.       
\end{abstract}

\section{Introduction}

The discovery of cubic structures with vivid photonic properties in a wide range of biological systems, such as those responsible for the   structural colouring seen in the wing scales of butterflies \cite{Saranathan_2010,Michielsen_2008,Gerd_2011,Wilts_2017} and in beetle cuticles \cite{Galusha_2008,Wilts_2012}, has driven a renewed interest in mathematics of \emph{minimal surfaces} -- surfaces which, everywhere, have zero mean curvature. The optical properties in these structures comes from the periodic nature of their bulk structure. Likewise, minimal surfaces that have three independent translational symmetries are known as \emph{triply-periodic minimal surfaces} (TPMS). 
Interfaces found in a plethora of soft matter systems --  including diblock copolymers \cite{Matsen_1994}, lyotropic liquid crystals, micellar solutions and lipid membranes \cite{Squires_2013} -- form structures approximating the geometry of three common triply-periodic minimal surfaces: the Schwarz D (aka the diamond surface or D-surface), the Schwarz P (aka the primitive or P-surface) \cite{Schwarz_1890} and the gyroid \cite{Schoen_1970}. Collectively, these structures are called cubic mesophases \cite{Matsen_1994,Gerd_2007,Luzzati_1967,Scriven_1976}.

The study of minimal surfaces has a rich and storied history, originating with Joseph-Louis Lagrange (1736-1813) \cite{Lagrange_1760}. He used the calculus of variations to find surfaces of minimum area with a given boundary curve. 
Based on this notion a minimal surface is defined as a critical point of the area functional.
This is equivalent to the condition that the mean curvature be everywhere zero on the surface. The Belgian physicist, Joseph Antoine Ferdinand Plateau (1801-1883) generated and studied several minimal surfaces with carefully designed experiments using soap films \cite{Plateau_1873}. 
The problem of proving the existence of a minimal surface with a given boundary curve is known as Plateau's problem. The Plateau problem for simple closed boundary curves in dimensions greater than or equal to two was solved independently by Jesse Douglas \cite{Douglas_1931} and Tibor Rad\'{o} \cite{Rado_1930}. Tibor Rad\'{o} also proved the uniqueness of a solution to Plateau's problem for all the boundary curves that can be mapped bijectively to a planar convex curve by either parallel or central projection \cite{Rado_1932}.

A triply-periodic minimal surface is a minimal immersion in $\mathbb{R}^3$ with a crystallographic or space group symmetry \cite{Paufler_2007}.  
The very first intersection-free (or embedded) triply-periodic minimal surface, the diamond surface was discovered by H. A. Schwarz \cite{Schwarz_1890}. It is named after the crystallographic structure of diamond as they share same set of symmetries and topology. Schwarz developed three more examples of intersection-free triply-periodic minimal surfaces, namely Schwarz's Primitive surface (or P surface), Schwarz's Hexagonal surface (or H surface) and the CLP surface \cite{Schwarz_1890}. Later his student E. R. Neovius discovered the C(P) surface also known as Neovius' surface \cite{Neovius_1883}. 

In the late 1960's, Alan Schoen explored the concept of labyrinth-spanning networks of a triply-periodic minimal surface and coined the term skeletal graph for such networks. 
Triply-periodic minimal surfaces are bicontinuous, meaning that the surface divides space into two disjoint region.
He imagined simultaneously inflating tubular neighbourhoods around both skeletal graphs. He postulated that a minimal surface would occur when these two inflationary regions collided. 
Using this notion, he constructed plastic models of seventeen examples of intersection-free triply-periodic minimal surfaces of which twelve were unknown at that time \cite{Schoen_1970}. He found a parametrisation for the gyroid using the theory of the Weierstrass-Enneper representation  \cite{Weierstrass_1903,Enneper_1868,Nitsche_1989} and Bonnet transformation \cite{Bonnet_1853}. Herman Karcher later proved the existence of the rest of Schoen's triply-periodic minimal surfaces  \cite{Karcher_1989}. In 1996, Karsten Gro\ss e-Brauckmann $\&$ Meinhard Wohlgemuth rigorously established that the gyroid is a minimal embedding \cite{KGB_i_1996}\footnote{They also proved embedding of the lidinoid, a hexagonal triply-periodic minimal surface, which, like the gyroid lacks straight lines and mirror planes.}. Recent additions to the list of triply-periodic minimal surfaces are the triply-periodic Costa surfaces \cite{Batista_2003} and three single-parameter families of triply-periodic minimal surfaces obtained by deforming the gyroid (tetragonal and rhombohedral deformation) and the lidinoid (rhombohedral deformation) \cite{Weyhaupt_2008}.

Over the years, researchers have proposed numerous ways to find new triply-periodic minimal surfaces: approaches based on space group or crystallographic symmetries by Fischer $\&$ Koch \cite{Koch_1988,Koch_i_1989,Koch_ii_1989,Fischer_i_1989,Fischer_ii_1989} and Lord $\&$ Mackay \cite{Lord_2003}, a method based on a Schwarz-Christoffel formula for periodic polygons in the plane by Fujimori $\&$ Weber \cite{Fujimori_2009}, an approach based on tilings of hyperbolic plane by Sadoc $\&$ Charvolin \cite{Sadoc_1989}, an algorithm based on the conjugate surface method and the concept of discrete minimal surfaces by Pinkall, Polthier $\&$ Karcher \cite{Pinkall_1993,Karcher_1996}, a method based on Schwarz triangular tilings of the two sphere by Fogden $\&$ Hyde \cite{Fogden_i_1992,Fogden_ii_1992}. Gandy et al. construct explicit parametrisations 
 \cite{Gandy_1999, Gandy_i_2000, Gandy_ii_2000}; such analytical representations are highly accurate compared with periodic nodal surface approximations  \cite{Gandy_2001}. Kenneth Brakke developed Surface Svolver \cite{Brakke_1992}, a software package uses a conjugate gradient solver to minimise energy functionals on a surface. Using crystallographic symmetries and boundary conditions,  Surface Evolver can numerically generate intersection-free translational unit cells of triply-periodic minimal surfaces \cite{KGB_1997, Hyde_2009}.

One of the most successful methods constructs triply-periodic minimal surface using a bottom-up approach by posing it as a Plateau problem for an appropriate boundary curve \cite{Karcher_1996}. If the boundary curve consists of \emph{in-surface symmetries} namely straight line segments, which give rise to two-fold rotational symmetries in the surface, and curves lying on mirror planes then the minimal surface patch bound within such a curve can be rotated and reflected to generate the translational unit of the triply-periodic minimal surface. 
The above method fails if a triply-periodic minimal surface has no in-surface symmetries. This is because unlike surfaces with in-surface symmetries there are no surface patches bounded by a special set of curves i.e. straight line segments and curves lying on planes of reflection symmetry. 

Here, we tackle the problem of constructing triply-periodic minimal surfaces with no in-surface symmetries. Our algorithm is a top-down approach, combining a numerical construction in Surface Evolver (based on a procedure suggested by de Campo, Hyde et al \cite{Hyde_2009} in the context of polycontinuous phases.) and the complex analysis of minimal surfaces using Weierstrass-Enneper representation. We use this method to describe and parametrise a new one-parameter family of chiral triply-periodic minimal surfaces with hexagonal symmetry. We call our new family of surfaces the QTZ-QZD family, named for their parent networks, the quartz or $\bf{qtz}$ network and its proper dual, the $\bf{qzd}$ network.

Our analytical representation of the QTZ-QZD family of surfaces is based on a method developed by Fogden and Hyde in a series of papers from 1992 \cite{Fogden_i_1992,Fogden_ii_1992}. In their first paper \cite{Fogden_i_1992}, they introduced the notion of regular and irregular triply-periodic minimal surfaces -- a classification scheme based on the behavior of the Gauss map in the vicinity of \emph{flat points} -- the points of zero gaussian curvature. Based on this scheme they unified several triply-periodic minimal surfaces under the regular class establishing an elegant connection between Schwarz triangular tilings of the unit sphere and the existence of a triply-periodic minimal surface with a given crystallographic symmetry. They focused only on the regular class in first two papers \cite{Fogden_i_1992,Fogden_ii_1992}, later in 1993 Fogden extended the algorithm to parametrise several irregular surfaces \cite{Fogden_iii_1993,Fogden_2_1994,Fogden_1_1994}. Nevertheless all of these triply-periodic minimal surfaces have in-surface symmetries. Our method is solely based on a pair of dual space graphs.
 Therefore, it can be used to generate both regular and irregular classes of surfaces regardless of their in-surface symmetries. Here, we highlight this by explicitly constructing the QTZ-QZD family of surfaces which is irregular and has no in-surface symmetries.

The QTZ-QZD family of surfaces are uniquely (up to isometries of $\mathbb{R}^3$) determined by a continuous parameter: the ratio of the chiral pitch to the hexagonal translational period.
In the following section, we give a brief outline of our algorithm for parametrising triply-periodic minimal surfaces. In section \ref{sec:tpms}, we introduce the Weierstrass-Enneper representation of minimal surfaces. A major part of this section is dedicated to deriving relations that govern different entities in the algebraic equation satisfied by the Weierstrass function; this is done by analyzing the geometric and topological properties of the surface. In section \ref{sec:fps}, we compute the branch points of the Weierstrass function up to two unknown parameters. In section \ref{sec:wf}, firstly, we calculate the Weierstrass function for the QTZ-QZD family of surfaces, and secondly, we solve the period problem to completely determine the family.


\begin{figure}[ht!]
\centering
\begin{subfigure}{.5\textwidth}
\labellist
\small\hair 2pt
\pinlabel $\bf{(a)}$ at 100 700
\endlabellist
  \centering
  \includegraphics[width=.95\linewidth]{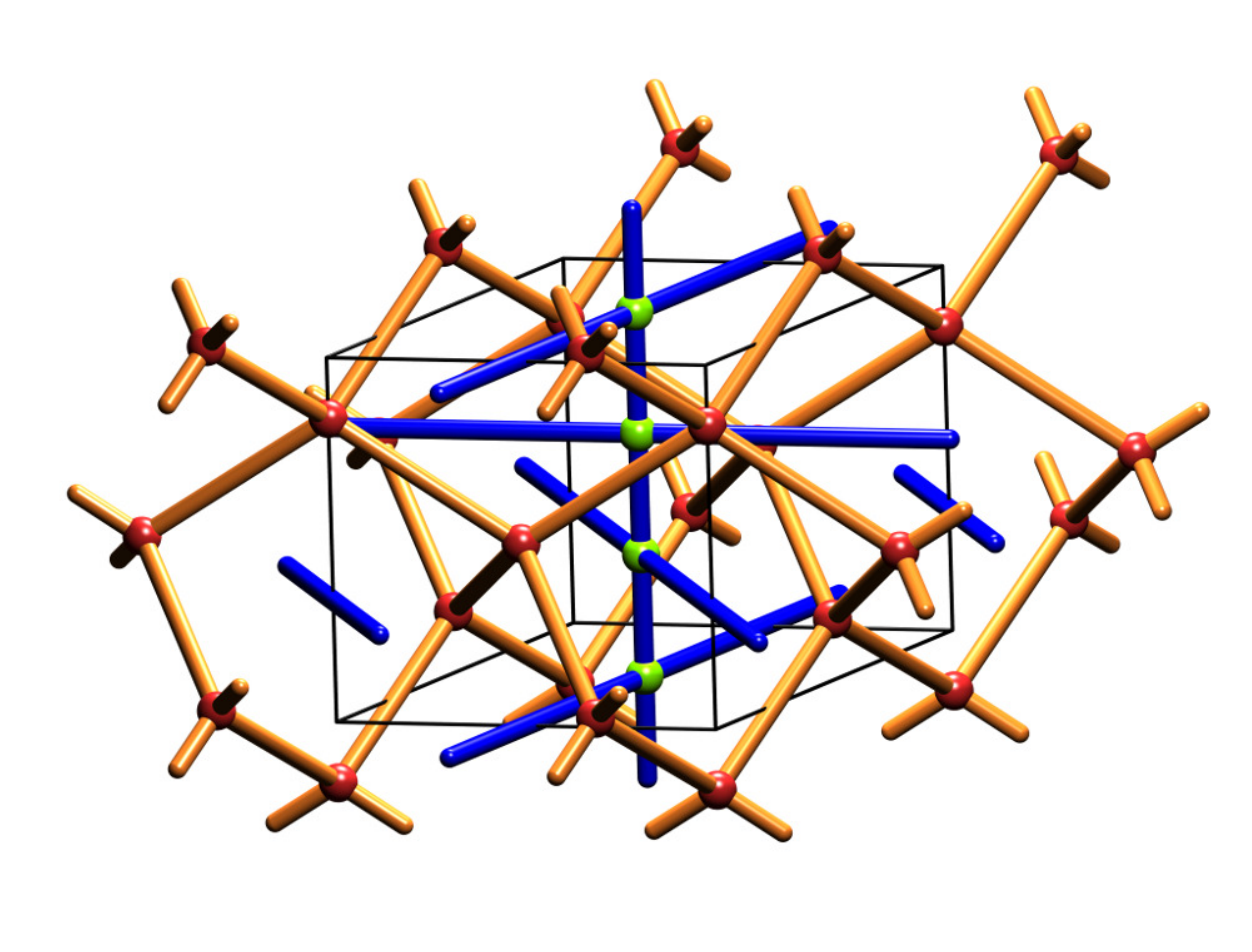}
  \label{fig:sub1}
\end{subfigure}%
\begin{subfigure}{.5\textwidth}
\labellist
\small\hair 2pt
\pinlabel $\bf{(b)}$ at 100 700
\endlabellist
  \centering
  \includegraphics[width=.95\linewidth]{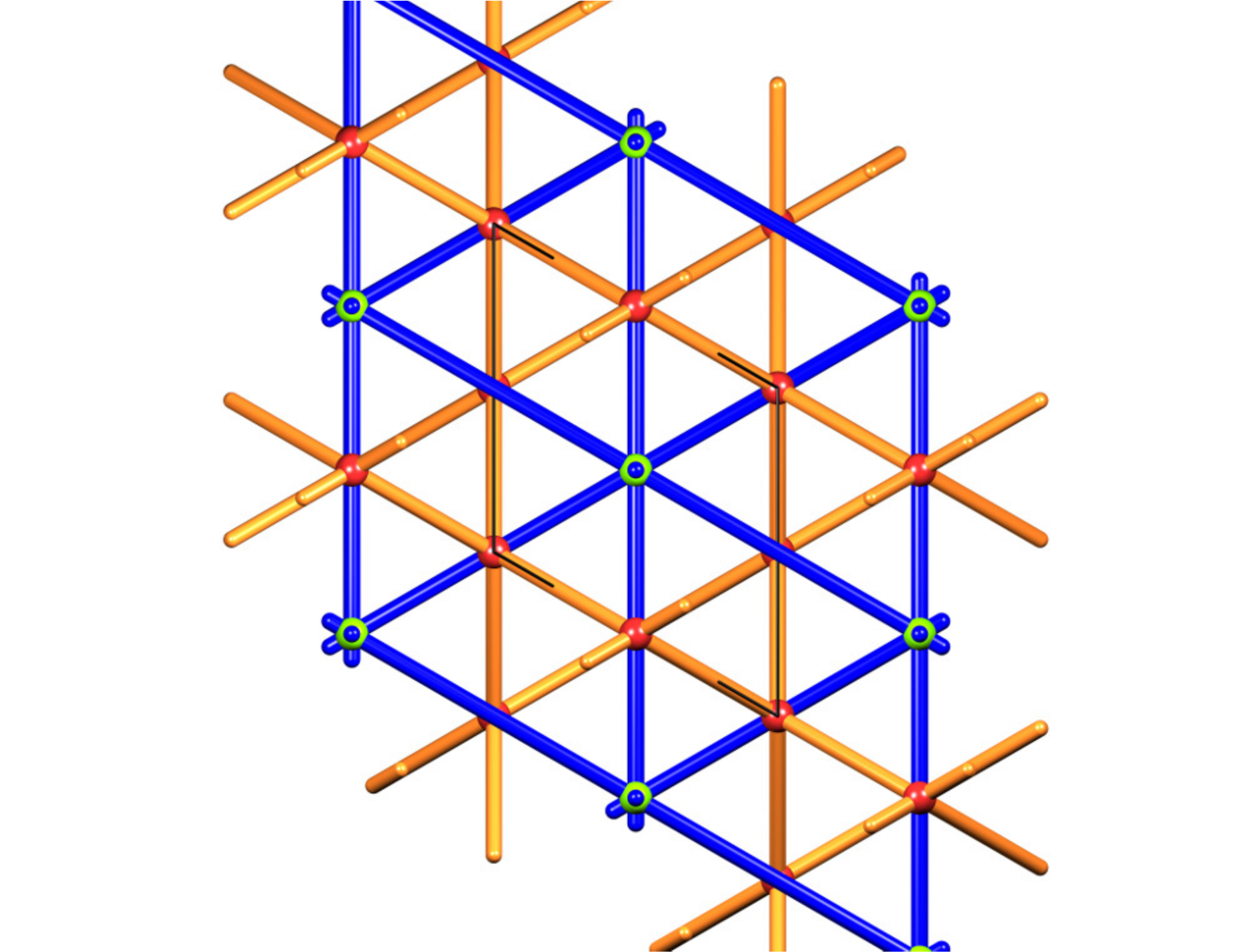}
  \label{fig:sub2}
\end{subfigure}
\caption{The pair of interpenetrating dual 3D graphs that lead to the QTZ-QZD surface family, the quartz (qtz) network and its proper dual {qzd} network (blue). (a) A top view; the translational unit cell is shown in black, (b) a side view; the translational unit cell of the space group lattice $P6_222$ (black).}
\label{fig:test}
\end{figure}

\section{An overview of the algorithm used to parametrise the QTZ-QZD family of surfaces}
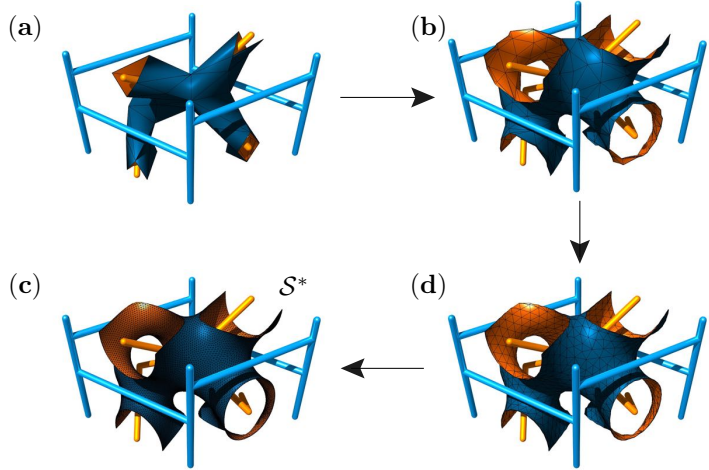
\begin{figure}[ht!]
\labellist
\small\hair 2pt
\pinlabel $\bf{(c)}$ at 10 200
\pinlabel $\bf{(d)}$ at 400 200
\pinlabel $\bf{(b)}$ at 400 450
\pinlabel $\bf{(a)}$ at 10 450
\pinlabel $\mathcal{S}^*$ at 270 200
\endlabellist
\centering
\includegraphics[scale=0.4]{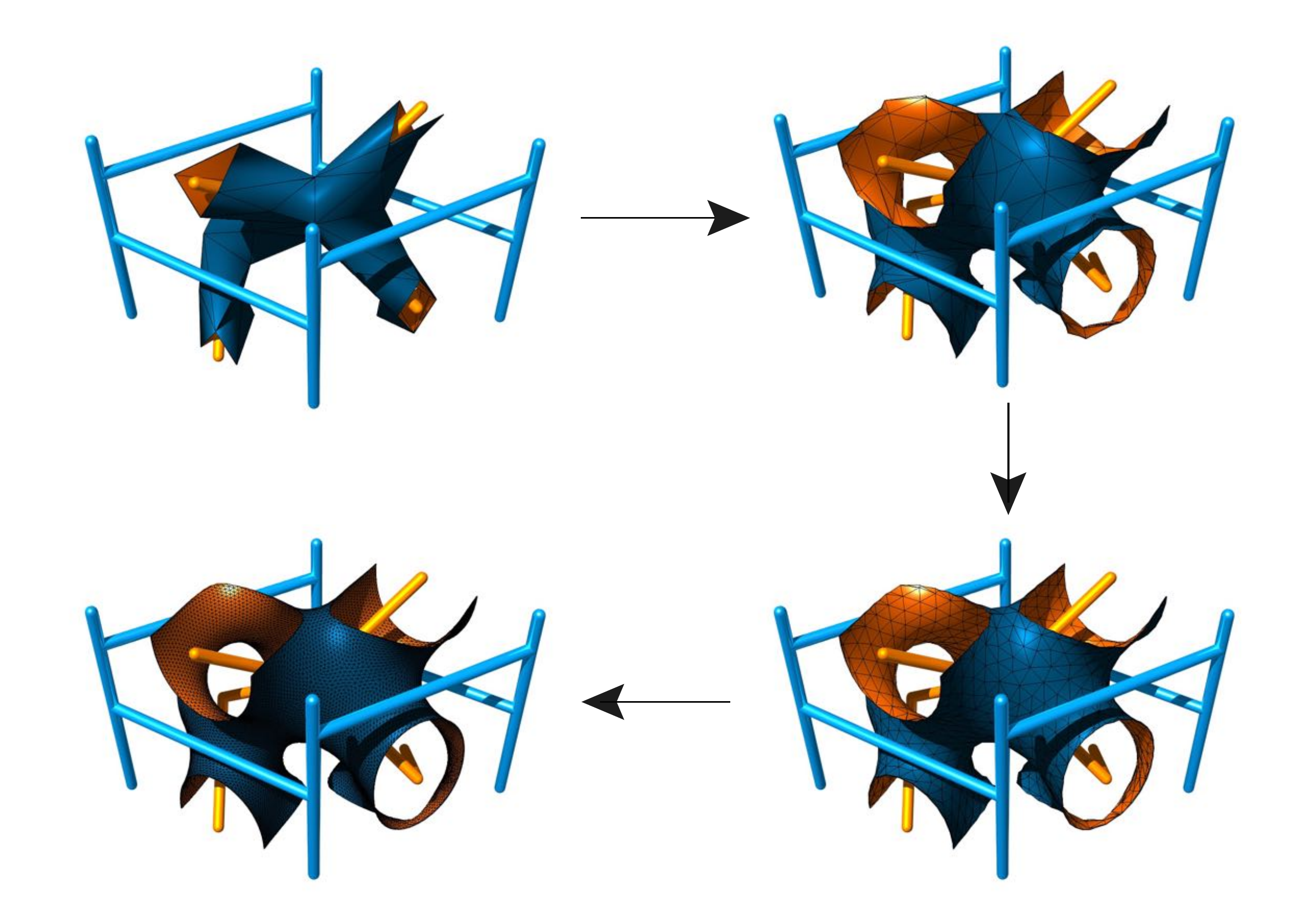}
\caption{Four different stages illustrating the minimisation of the area functional using conjugate gradient descent method in Surface Evolver \cite{Brakke_1992}. See Appendix A2.2. The boundary conditions respect the translational periodicities imposed by the $P6_222$ space group symmetry with lattice parameters $c= a=1$. The blue coloured network is the {qzd}  network and the orange coloured network is the {qtz} (quartz) network. Figure (a) shows a triangulated tubular neighbourhood of the quartz network and figure (c) shows the final surface, $\mathcal{S}^*$ after convergence.}
\label{fig:test1}
\end{figure}

Motivated by Alan Schoen's method of generating infinitely periodic minimal surfaces in \cite{Schoen_1970}, we aim to find a unique triply-periodic minimal surface starting with a pair of dual interpenetrating 3D networks i.e. the qtz network \cite{OKeeffe_2008,DelgadoF_2002,DelgadoF_1_2003,DelgadoF_2_2003} and the qzd network\footnote{These and many other three-dimensional networks and their dual networks are documented in RCSR \cite{OKeeffe_2008}.} \cite{OKeeffe_2008,DelgadoF_2_2003,Rosi_2005} (shown in Figure \ref{fig:test}). Our goal is to construct an intersection-free triply-periodic minimal surface that divides the ambient space into two disjoint multiply connected three-dimensional domains such that one of them is homeomorphic to a tubular neighbourhood of the qtz network whilst the other is homeomorphic to a tubular neighbourhood of qzd network; such a partition gives rise to the notion of labyrinth-spanning networks called the skeletal graphs \cite{Schoen_1970}. We use a method put forth by Karsten Gro\ss e-Brauckmann in \cite{KGB_1997} to obtain an intersection-free triangulated surface. Using a triangulated tubular neighbourhood of the translational unit cell of the qtz network and the method of conjugate gradient descent, we get a triply-periodic discrete surface; four stages in the evolution of the numerical algorithm implemented in Surface Evolver are shown in Figure \ref{fig:test1}. It follows from the construction protocol that this surface is an approximation to a minimal surface embedded in a three torus. For a more elaborate discussion on the numerical aspects of the algorithm refer to the Appendix A2.

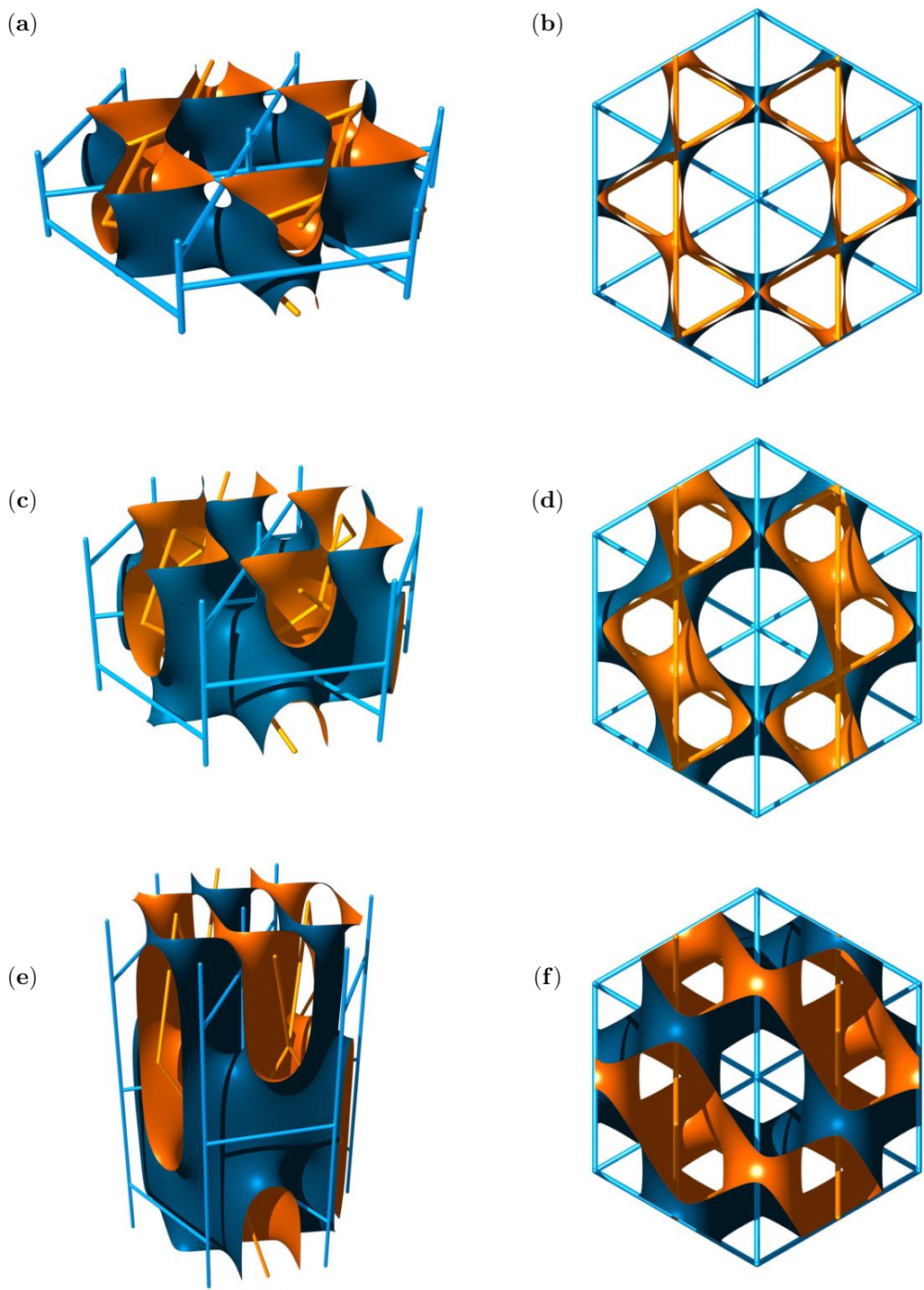
\begin{figure}
\labellist
\small\hair 2pt
\pinlabel $\bf{(a)}$ at 10 1100
\pinlabel $\bf{(b)}$ at 450 1100
\pinlabel $\bf{(c)}$ at 10 700
\pinlabel $\bf{(d)}$ at 450 700
\pinlabel $\bf{(e)}$ at 10 300
\pinlabel $\bf{(f)}$ at 450 300
\endlabellist
\centering
\includegraphics[scale=.5]{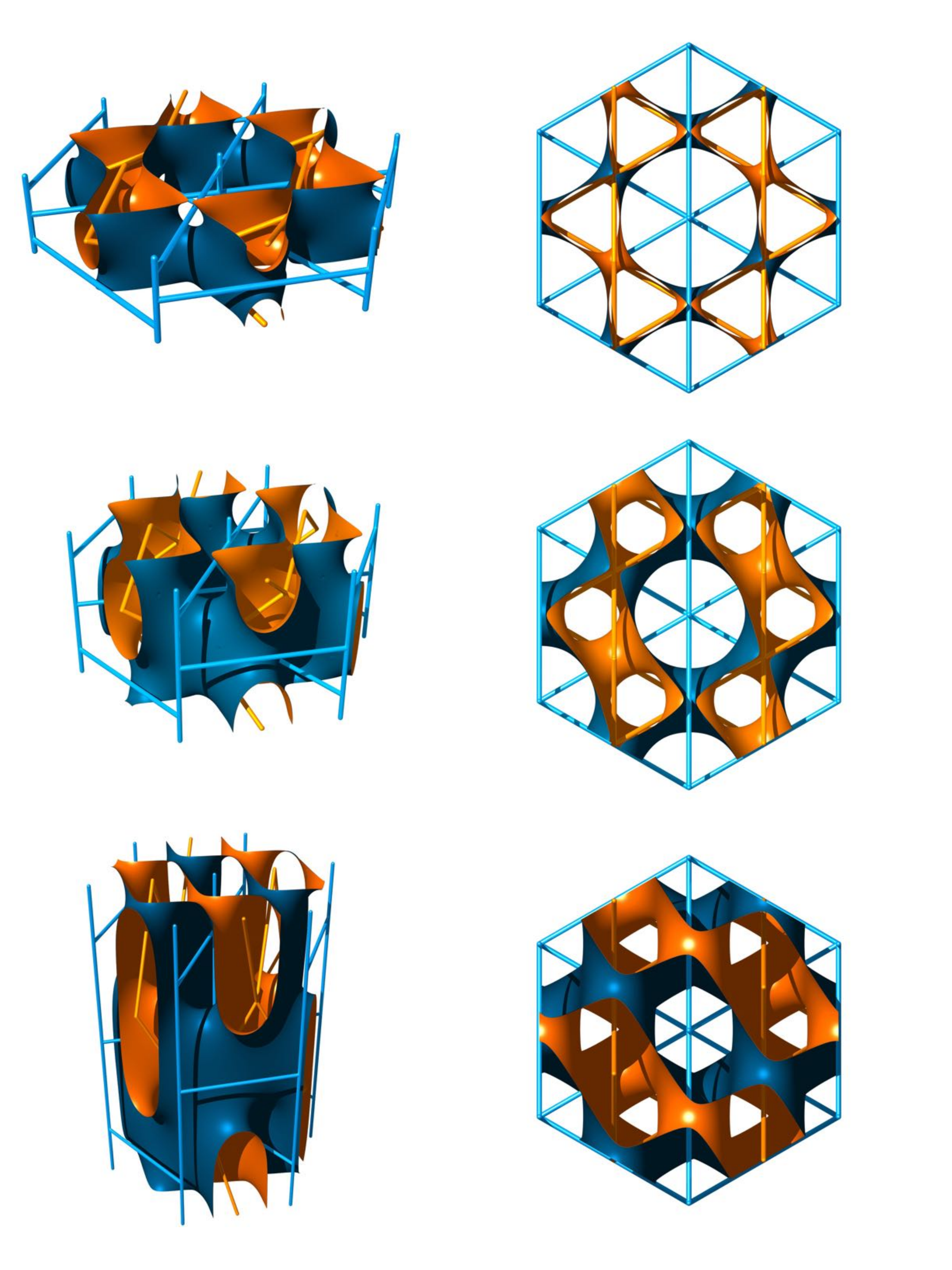}
\caption{Three members of the QTZ-QZD TPMS family with different values of $\rho$. (a) A side view, (b) a top view for $\rho\approx$ 0.445. (c) A side view, (d) a top view for $\rho\approx$ 0.994. (e) A side view, (f) a top view for $\rho\approx$ 2.887. The labyrinth enclosed within inner side (orange coloured) of the surface is spanned by the qtz network.}
\label{fig:test8}
\end{figure}   

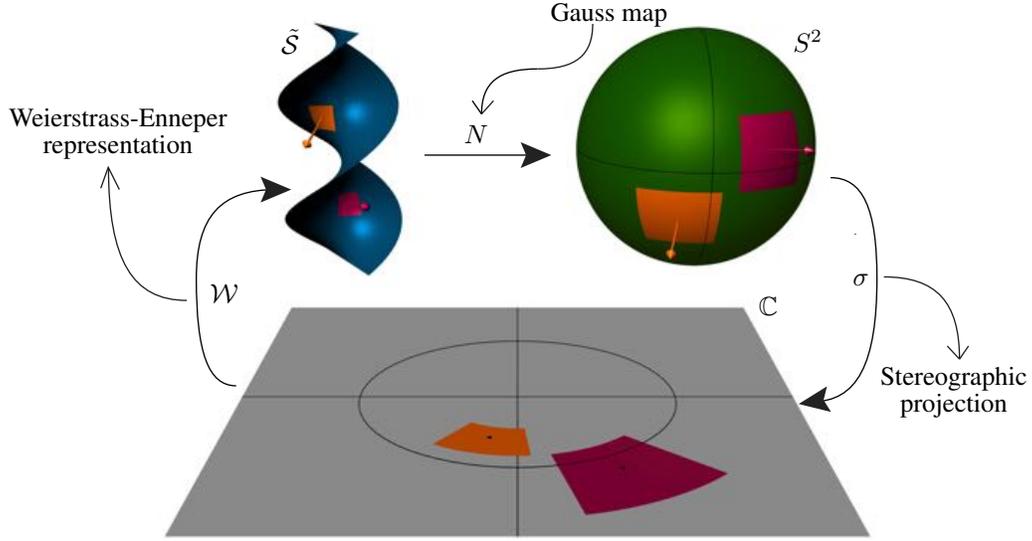
\begin{figure}[ht!]
\labellist
\small\hair 2pt
\pinlabel $\mathcal{W}$ at 110 190
\pinlabel $N$ at 300 310
\pinlabel $\sigma$ at 590 200
\pinlabel $\mathcal{\tilde{S}}$ at 160 380
\pinlabel $S^2$ at 550 380
\pinlabel $\mathbb{C}$ at 520 180
\pinlabel $\textrm{Gauss map}$ at 400 400
\pinlabel $\textrm{Weierstrass-Enneper}$ at 30 320
\pinlabel $\textrm{representation}$ at 30 300
\pinlabel $\textrm{Stereographic}$ at 660 125
\pinlabel $\textrm{projection}$ at 660 105
\endlabellist
\centering
\includegraphics[scale=0.5]{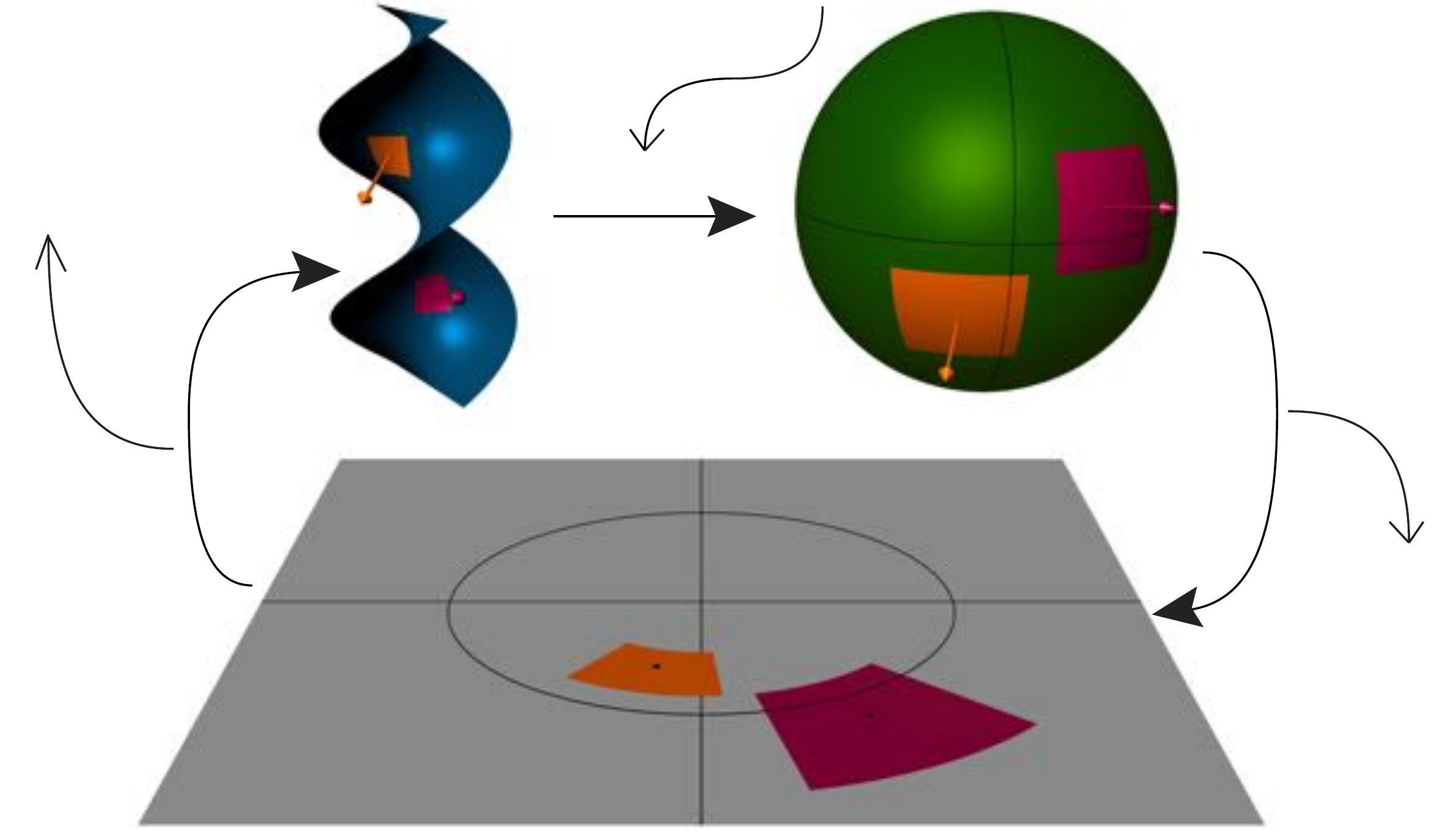}
\caption{A schematic diagram illustrating the key steps involved in parametrising a minimal surface. The Weierstrass-Enneper representation of the minimal surface $\mathcal{\tilde{S}}$ is obtained by solving the period problem based on the symmetries, Gauss map and curvature data of $\mathcal{S}^*$. Note that in our case $\mathcal{\tilde{S}}$ refers to a surface from the QTZ-QZD surface family and the minimal surface in the figure, helicoid is used just to represent a generic minimal surface.}
\label{fig:test2}
\end{figure}

The numerical construction by itself is not sufficient to establish that the discretised surface corresponds to a minimal surface. To establish that the discretised surface corresponds to a embedded minimal surface with rigor we can proceed in two ways: \textrm{(i)} we can establish by numerical analysis that the limiting surface under the limit of mesh size going to zero is a minimal surface, \textrm{(ii)} we can find a parametrisation of the corresponding triply-periodic minimal surface using techniques from differential geometry and the complex analysis of minimal surfaces. We take the latter route and prove the existence of the QTZ-QZD family of minimal surfaces by explicit construction. Three members of which are shown in \ref{fig:test8}.  A schematic of our approach is illustrated in Figure \ref{fig:test2}. The method developed here heavily relies on the Weierstrass-Enneper representation proved by Alfred Enneper and Karl Weierstrass back in 1860s  \cite{Weierstrass_1903,Enneper_1868}. The Weierstrass-Enneper representation has been a critical tool in discovering new intersection-free triply-periodic minimal surfaces \cite{Schoen_1970, Fogden_i_1992,Fogden_ii_1992,Fogden_iii_1993,Lidin_1987,Lidin_1990,Fogden_1_1994,Fogden_2_1994,Batista_2003,Weyhaupt_2008}. In the rest of the paper, we describe various components of the Weierstrass-Enneper representation and how the numerical data extracted from the discretised surface leads to an embedded triply-periodic minimal surface consistent with the topology and symmetries of the qtz/qzd networks.   


\section{Triply-periodic minimal surfaces}\label{sec:tpms}
Our interest in studying triply-periodic minimal surfaces is rooted in their physical manifestations and engineering applications so we restrict ourselves only to the \emph{properly embedded} surfaces. A surface $\mathcal{M}:\mathcal{D}\rightarrow X \subset \mathbb{R}^3$ is said to be \emph{proper} when the preimage of compact sets are compact. Such triply-periodic minimal surfaces are devoid of self-intersections and they divide the ambient space into two disjoint connected interpenetrating three-dimensional labyrinthine domains. Properly embedded surfaces avoid unphysical pathological situations such as infinite labyrinths in a finite volume in $\mathbb{R}^3$. 

Let us denote the triply-periodic minimal surface by $\mathcal{M}$. The set of isometries of $\mathbb{R}^3$ with standard euclidean metric preserving $\mathcal{M}$ forms a space group. Let us denote the subgroup of translations by $\Lambda$. The surface patch enclosed inside the translational unit cell of the space group denoted by $\mathbb{R}^3/\Lambda$ is the quotient of $\mathcal{M}$ with respect to $\Lambda$. Therefore the group action of $\Lambda$ on $\mathcal{M}$ gives us the smallest periodic domain of $\mathcal{M}$ denoted by $\mathcal{M}_{\Lambda}$, where $\mathcal{M}_{\Lambda}\equiv \mathcal{M}/\Lambda\subset\mathbb{R}^3/\Lambda\equiv T^3.$
$T^3$ denotes a three-torus equivalent to the translational unit cell of the space group. The boundary components of $\mathcal{M}_{\Lambda}$ are defined as the curves the intersection of $\mathcal{M}$ with the boundary of $\mathbb{R}^3/\Lambda\equiv T^3$. The generators of $\Lambda$ uniquely identify all the boundary components of $\mathcal{M}_{\Lambda}$ in pairs. Hence, $\mathcal{M}_{\Lambda}$ embedded in $T^3$ is compact and has a finite total gaussian curvature. According to Theorem 1.1 in Appendix A, there exists a conformal parametrisation\footnote{A heuristic derivation of \ref{eqn:wer} is given in the Appendix I(a) in the ESM.}  for $\mathcal{M}/\Lambda$\footnote{From here onwards, we assert that 
all functions of the surface here discussed, including the coordinate functions and Gauss map restricted to $\mathcal{M}_{\Lambda}$,
are well defined and they can be extended to the lift of $\mathcal{M}_{\Lambda}$ in $T^3$ to $\mathcal{M}$ in $\mathbb{R}^3$ giving the Gauss map and the parametrisation of $\mathcal{M}$ \cite{Hoffman_1990}.} given by
\begin{align}\label{eqn:wer}
\mathcal{W} & : \mathcal{D}\rightarrow\mathcal{M}/\Lambda\subset\mathbb{R}^3/\Lambda\nonumber\\
 x_1(\xi) & = \textrm{Re}\left(e^{i\theta}\int_{\xi}\,(1-z^2)\,R(z)\,dz\right) + c_1\nonumber\\ 
 x_2(\xi) & = \textrm{Re}\left(e^{i\theta}\int_{\xi} \,i(1+z^2)\,R(z)\,dz\right) + c_2\nonumber\\
 x_3(\xi) & = \textrm{Re}\left(e^{i\theta}\int_{\xi}\,2z \,R(z)\,dz\right)+ c_3. 
\end{align}
where the function $R$ is the Weierstrass function, $\theta$ is the Bonnet angle \cite{Bonnet_1853} and $\mathcal{D}\subset\mathbb{C}\cup\{\infty\}$ is simply connected. Above parametrisation gives a family of isometric minimal surfaces known as the \emph{associate family}. The multiplication by a complex number of unit magnitude is known as the Bonnet transformation. The Bonnet transformation is an isometry and the continuous parameter $\theta$ is called the Bonnet angle. Isometric minimal surfaces differing by a Bonnet angle of $\pi/2$ are said to be conjugate or adjoint to each other. Some examples of adjoint minimal surfaces are Catenoid, Helicoid and Schwarz's P surface, Schwarz's D surface. The Bonnet transformation has proven to be a valuable tool in generating new triply-periodic minimal surfaces \cite{Pinkall_1993,Karcher_1996}. The gyroid \cite{Schoen_1970} and the lidinoid \cite{Lidin_1990} are generated by Bonnet transformation of Schwarz's D surface and Schwarz's H surface respectively.


The Gauss map $N:\mathcal{M}_{\Lambda}\rightarrow S^2$ is a function from a surface to the unit sphere. It maps a point on the surface to the unit normal vector at that point.
The number of sheets $s$ in the Gauss map, branched cover of $S^2$, can be related to the genus $g$ of $\mathcal{M}_{\Lambda}$ and its gaussian curvature $K$ by
\begin{align*}
\oint_{\mathcal{M}_{\Lambda}} KdA & = -s\oint_{S^2} dA\\
2\pi\chi & =  -4\pi s \\
g -1 & = s 
\end{align*}
where $K$ is the gaussian curvature
 and $\chi$ is the Euler characteristic of $\mathcal{M}_{\Lambda}$. The second equality above results from the \emph{Gauss-Bonnet theorem} \cite{Gauss_1827,Bonnet_1848,vonDyck_1888}. Thus by composing the Gauss map and the stereographic projection we get a $(g-1)$ sheeted branched cover of $\mathbb{C}\cup\{\infty\}$. In fact $\mathcal{M}_{\Lambda}$ parametrised by \ref{eqn:wer} is an immersion of a genus $g$ Riemann surface into $T^3$. The Riemann surface of interest is completely determined by the pair of functions $(\xi,R(\xi))$ satisfying an algebraic equation \cite{Springer_1957} (details in the following section), so the domain of integration $\mathcal{D}$ in \ref{eqn:wer} is a subset of this Riemann surface.


\subsection{The Weierstrass function}
The Weierstrass function satisfies the following algebraic equation for every $\xi\in\mathbb{C}\cup\{\infty\}$
\begin{equation}\label{eqn:pole}
\sum\limits_{m=0}^s a_m(\xi)R^m = 0
\end{equation}
where $a_m(\xi)$ is a polynomial over $\mathbb{C}$. Since the image of the Gauss map is a $s$-sheeted branched cover of $S^2$ the degree of this equation must be equal to $s$. 
The coefficient polynomials can be expressed as
\begin{equation}\label{eqn:pol_coef}
a_m(\xi) = \alpha_mP_m(\xi)\prod_{i=1}^n (\xi-\xi_i)^{q_{m,i}}
\end{equation}
where $q_{m,i}\geq 0$ is an integer, $P_m(\xi_i)$ is a polynomial, $\alpha_m$ is a complex constant for all $m\in\{0,1,2,...,s\}$ and $\{\xi_i\}_{i=1}^n$ is the set of branch points. 
The equation (1.12) in Appendix I a(ii) in the ESM gives an expression of the gaussian curvature in terms of the Weierstrass function.
\begin{equation}\label{eqn:gc}
K(\xi) = -\frac{4}{(1+|\xi|^2)^4|R(\xi)|^2}
\end{equation}
Using \ref{eqn:gc} we can infer a couple of things about the polynomial coefficients in \ref{eqn:pole}: \textrm{(i)} because of the Gauss-Bonnet theorem the total integrated gaussian curvature is finite and since it is non-positive through out the surface it must be finite everywhere, which implies $R(\xi)$ is non-zero for all $\xi\in\mathbb{C}\cup\{\infty\}$. Since $a_0(\xi)$ cannot vanish, $a_0$ must be independent of $\xi$. Therefore $a_0(\xi)=\alpha_0$,
where $\alpha_0$ is a complex valued constant. \textrm{(ii)} The singularities of the Weierstrass function are its branch points. Whenever $R(\xi)$ diverges the gaussian curvature vanishes, so the product of the roots of \ref{eqn:pole} equal to the ratio of the polynomials $a_s(\xi)$ and $a_0(\xi)$ must vanish at the branch points which implies that $a_s(\xi)$ must vanish at the branch points. Without loss of generality we assert $P_s(\xi)=1,$ which implies that $a_s(\xi) = \alpha_s\prod_{i=1}^n (\xi-\xi_i)^{q_{s,i}}.$

\subsubsection{The limiting behaviour at the branch points}
The \emph{Riemann-Hurwitz formula} \cite{Griffiths_1989} gives a relation between the total branch point order $W$ of the Weierstrass function, the number of sheets $s$ and the genus $g$ of the Riemann surface.
\begin{equation}
W  = 4(g-1) = 4s. 
\end{equation}
In general, a branch point in the complex plane, $\mathbb{C}$ can have more than one branched local neighbourhoods, distinct or non-distinct (with respect to the number of sheets pinned at the branch point). This way every branch point $\xi_i\in\mathbb{C}$ can be associated with a set $\mathcal{B}_i  = \{b_{ij},N_{ij}\}_{j=0}^{s-1}$
subject to the constraints
\begin{align}\label{eqn:bps}
&\sum\limits_{j=0}^{s-1}N_{ij}(b_{ij}+1)  = s\nonumber\\
&\sum\limits_{i=1}^n\sum\limits_{j=0}^{s-1}N_{ij}b_{ij}  = W = 4s.
\end{align}

The set $\mathcal{B}_i$ gives the \emph{branch point structure} in the neighbourhood of the $i^{th}$ branch point, $\xi_i\in\mathbb{C}$. When viewed as part of the branched cover there are $b_{ij}+1$ sheets pinned at the $j$th copy of the $i^{th}$ branch point and there are $N_{ij}$ such copies amounting to a total of $\sum_{j}N_{ij}(b_{ij}+1)$ sheets over $\xi_i\in\mathbb{C}$. A visualisation of the branched cover is shown in Figure 2 of the ESM. The approach followed here starting from \ref{eqn:pol_coef} in the last section was originally introduced and systematised in a series of two papers by Fogden $\&$ Hyde from 1992 \cite{Fogden_i_1992,Fogden_ii_1992}. They proposed an algorithm for obtaining the Weierstrass-Enneper representation which lead to a classification scheme based on the branch point structure of the underlying Riemann surface. They coined the term \emph{regular} triply-periodic minimal surfaces for surfaces with $\mathcal{M}_{\Lambda}$ homeomorphic to a Riemann surface with $|\mathcal{B}_i|$ equal to one for all $i\in\{1,2,...,n\}$, i.e. $N_i = \frac{s}{b_i+1}\implies \sum_{i=1}^n \frac{b_i}{b_i+1} = 4.$
Fogden extended this formalism to irregular\footnote{The triply-periodic minimal surfaces that do not belong regular class are categorised as irregular surfaces.} triply-periodic minimal surfaces in \cite{Fogden_iii_1993} by considering the branch point structure to be more generic as given by \ref{eqn:bps}.

The function $R(\xi)$ in the neighbourhood of the $j^{th}$ copy of the $i^{th}$ branch point can be written as
\begin{equation}\label{eqn:near_bp}
 R(\xi)\sim\gamma_{ij}(\xi-\xi_i)^{-b_{ij}/(b_{ij}+1)}
\end{equation}
where $b_{ij}$ is the branch point order \cite{Fogden_iii_1993}. In this notation $j=0$ corresponds to an unbranched sheet, implying that $b_{i0}$ equals zero. By substituting $R(\xi)$ from \ref{eqn:near_bp} into \ref{eqn:pole}  we get for the $m^{th}$ term
\begin{equation}
 a_m(\xi)R(\xi)^m\sim\alpha_mP_m(\xi_i)\gamma_{i0}^m(\xi-\xi_i)^{q_{m,i}}\prod_{k=1,k\neq i}^n(\xi_i-\xi_k)^{q_{m,k}}.
\end{equation}
If there are $N_{i0}$ unbranched sheets then $\gamma_{i0}$ must take $N_{i0}$ distinct values one on each unbranched sheet. This implies in the limit $\xi \rightarrow \xi_{i}$
\begin{align}\label{eqn:bp_anls}
&\alpha_{N_{i0}}P_{N_{i0}}(\xi_i)\gamma_{i0}^{N_{i0}} (\xi-\xi_i)^{q_{N_{i0},i}}\prod_{k=1,k\neq i}^n(\xi_i-\xi_k)^{q_{N_{i0},k}} \sim \alpha_0 \implies q_{N_{i0,i}}  = 0.\nonumber\\
&\alpha_mP_m(\xi_i)\gamma_{i0}^m(\xi-\xi_i)^{q_{m,i}}\prod_{k=1,k\neq i}^n(\xi_i-\xi_k)^{q_{m,k}}\sim O(1)\implies q_{m,i} \geq 0\nonumber\\
&\alpha_mP_m(\xi_i)\gamma_{i0}^m(\xi-\xi_i)^{q_{m,i}}\prod_{k=1,k\neq i}^n(\xi_i-\xi_k)^{q_{m,k}}\sim 0\implies q_{k,i} > 0
\end{align}
for $i\in\{1,2,...,n\}$,  $m\in\{1,2,...,N_{i0}-1\}$ and $k\in\{N_{i0}+1,...,s\}$. 

\subsubsection{The asymptotic behaviour}
It follows from \ref{eqn:gc} that $R(\xi)$ must satisfy
\begin{equation}\label{eq:asymptotic_behaviour}
 \lim_{\xi \to \infty}R(\xi)\sim \gamma_{\infty}|\xi|^{-4}.
\end{equation}
It is always possible to ascertain that the function is finite at complex infinity. Therefore, in order to ensure the existence of $s$ distinct values of $\gamma_{\infty}$ in the limit $\xi\rightarrow\infty$, we must have 
\begin{equation}\label{eqn:qs}
a_s(\xi)R^s(\xi) + \alpha_0\sim O(1)
\implies\sum\limits_i^n q_{s,i} = 4s.
\end{equation}
The rest of the terms can either vanish or be of order one in the limit $\xi\rightarrow\infty$, hence
\begin{equation}\label{eqn:asmp_anls}
\sum\limits_i^n q_{m,i} + \textrm{deg}(P_m)\leq 4m
\end{equation}
for all $m \in\{1,2,...,s-1\}$.

\begin{figure}[ht!]
\centering
\includegraphics[scale=.2]{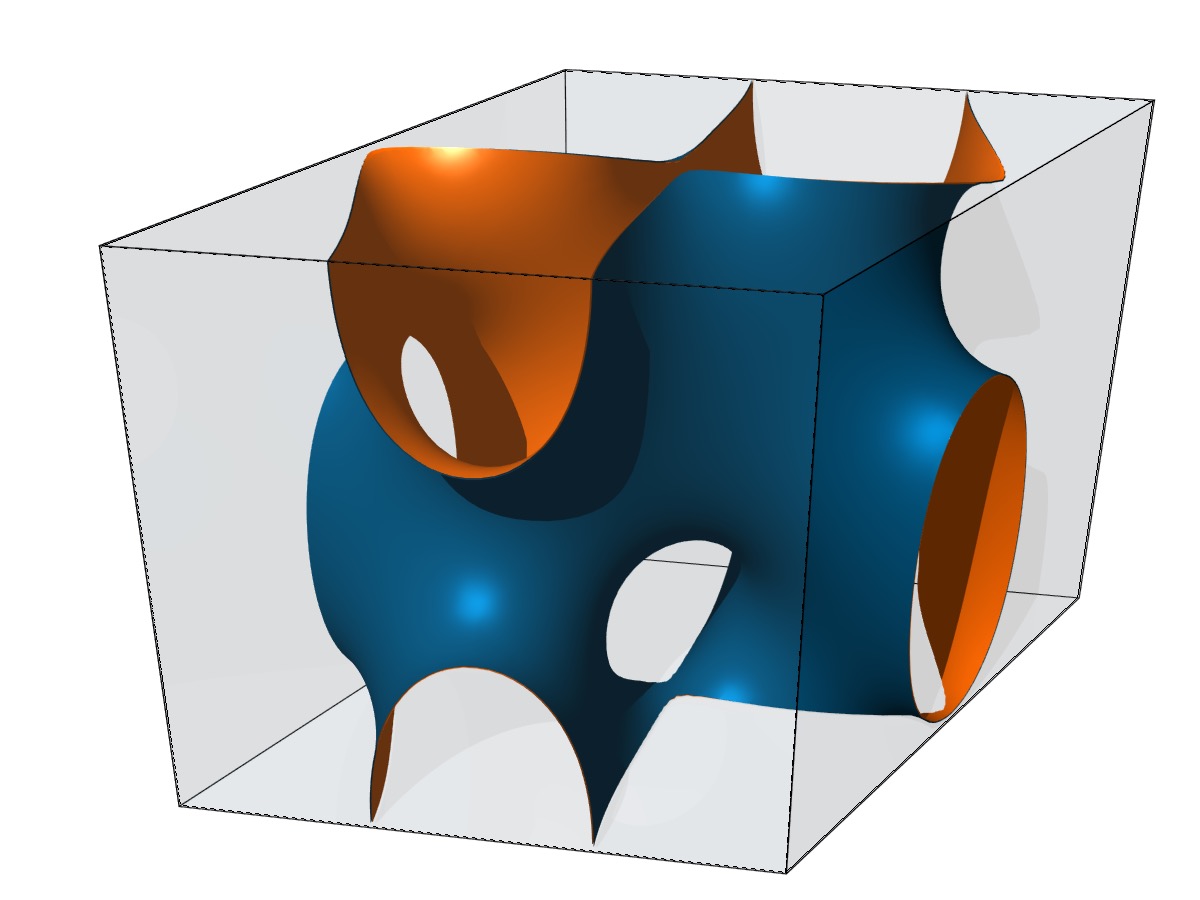}
\caption{The translational unit cell of $P6_222$ space group lattice of height $c$ and width $a$ such that $c/a\approx 1$, enclosed within is $\mathcal{M}_{\Lambda}$ for the QTZ-QZD surface with $\rho\approx1$.}
\label{fig:test3}
\end{figure}

\subsection{Space group of the QTZ-QZD family of surfaces}
The labyrinths of a triply-periodic minimal surface can be embedded with a unique pair of dual interpenetrating three dimensional networks called the \emph{skeletal graphs}. By definition, a triply-periodic minimal surface inherits the symmetries of its skeletal graphs \cite{Schoen_1970}. Therefore the space group preserving the QTZ-QZD family of surfaces is $P6_{2}22$. The $6_2$ represents a three-fold screw axis symmetry. The screw axis symmetry implies that there are no mirror symmetries (planar reflection symmetries) and that it is chiral. 

The primitive translational unit cell of an infinite crystallographic lattice can be generated using a smaller polyhedron with a minimal set of symmetry operations \cite{Paufler_2007}. A triply-periodic minimal surface endowed with a space group symmetry can be associated with a minimal surface patch (more on this in section \ref{sec:wf}\ref{sec:map}) enclosed within a minimal generating polyhedron. The minimal generating polyhedron for the space group lattice $P6_222$ shown in Figure \ref{fig:test8}(a) is a triangular slab of thickness $c/6$ and equilateral base of length $a$, where $c$ is the vertical period and $a$ is the horizontal period. The set of symmetry operations that act on the minimal surface patch to generate the translational unit $\mathcal{M}_{\Lambda}$ shown in Figure \ref{fig:test3} consists of: 
$(1)$ two two-fold rotational symmetries with nonintersecting axes parallel to the $xy$-plane and $(2)$ two two-fold rotational symmetries with axes parallel to the $z$-axis \cite{Dunbar_1988}. The action of these symmetries results in an effective three-fold screw axis parallel to the $z$-axis. 
 In the following, we summarise the procedure used to derive equations involving the parameters $\alpha_m$, $\textrm{deg}(P_{m,i})$ and $q_{m,i}$ by analyzing the action of the symmetry elements on the Weierstrass function.
    
\subsubsection{The action of rotations on the Gauss map $\&$ the Weierstrass function}
We choose a coordinate system with $z$-axis along the vertical direction. The invariance of $\mathcal{M}_{\Lambda}$ under the action of symmetry transformations implies mainly two things: \textrm{(i)} The coordinates of the minimal surface as given by the Weierstrass-Enneper representation \ref{eqn:wer} must be invariant, \textrm{(ii)} the polynomial equation \ref{eqn:pole} must be invariant, therefore the transformation of $a_m(\xi)$ for every value of m is completely determined by the transformation of $R(\xi)$.

The action of a rotation is determined by an orthogonal matrix. We multiply the Gauss map by this matrix to calculate the transformed Gauss map which in turn determines the transformation of the holomorphic Gauss map $\xi$. Then (i) leads us to the transformation of $R(\xi))$ and (ii) gives us the symmetries of the polynomial coefficients $a_m(\xi)$. In the following, we summarise the constraint equations obtained by acting with the two-fold rotational symmetry about an axis parallel to $z$-axis.
\begin{equation*}
\begin{split}
&\{\xi_i\} = \{-\xi_i\} \\
&\{\xi_{m,i}\} = \{-\xi_{m,i}\}
\end{split}
\quad\quad
\begin{split}
&R(\xi)  =  R(-\xi)\\
&a_m(\xi)  = a_m(-\xi)
\end{split}
\quad\quad
\begin{split}
(-1)^{\textrm{deg}(P_m) +\sum\limits_i^n q_{m,i}} = 1
\end{split}
\end{equation*}
Let $(\sin{\phi},\cos{\phi},0)$ be a two-fold rotation axis lying on the $xy$-plane. Then, 
\begin{eqnarray*}
&\displaystyle\xi^{4m} a_m\left(\frac{e^{-2i\phi}}{\xi}\right) = e^{-4i\phi m}a_m(\xi), \quad \xi^4 R(\xi) = e^{-4i\phi}R\left(\frac{e^{-2i\phi}}{\xi}\right)\\
&\displaystyle\{\xi_{m,i}\}  = \left\{\frac{e^{-2i\phi}}{\xi_{m,i}}\right\}, \quad \{\xi_i\} = \left\{\frac{e^{-2i\phi}}{\xi_i}\right\}\\
&\displaystyle(-1)^{\textrm{deg}(P_m)+\sum\limits_{i=1}^{n}q_{m,i}}e^{4i\phi m}\left(\prod_{i=1}^{n}\xi_i^{q_{m,i}}\right)\left(\prod_{j=1}^{\textrm{deg}(P_m)}\xi_{m,i}\right) = 1, \quad \textrm{deg}(P_m) +  \sum\limits_{i=1}^{n}q_{m,i} = 4m.
\end{eqnarray*}
The previous equations are derived in Appendix I(b) of the ESM.
 
\subsubsection{The screw axis symmetry}
There are in total three distinct screw axes in the translational unit (shown in Figure \ref{fig:test3}) among which two are along the labyrinth spanned by the {qtz} network and one is on the boundary of the unit cell along the labyrinth spanned by the {qzd} network. The symmetry operation involves rotating by $2\pi/3$ followed by a translation of one $c/3$, a third of the vertical period along the rotation axis (parallel to the $z$-axis). The transformation of the holomorphic Gauss map $\xi$ is given by $\xi^{\prime}=e^{i(2\pi/3)}\xi$
for an anti-clockwise screw axis. The pairing between the sense of rotation and the direction of translation is governed by the chirality of the skeletal graphs. The transformation of the $z$-coordinate from (3.1) gives $R(\xi) = e^{(4i\pi/3)}R(e^{(2i\pi/3)}\xi).$
%

\section{The flat points of the QTZ-QZD surface}\label{sec:fps}
\subsection{The branch point structure}
We find the genus of $\mathcal{M}_{\Lambda}$ to be four both using the skeletal graphs and by applying Gauss-Bonnet theorem to the discretised surface. Therefore, the number of sheets in the branched cover of the sphere is three. By substituting $s=3$ and $W=12$ in \ref{eqn:bps}, we get 
\begin{eqnarray}
 \sum\limits_{j=0}^2 N_{ij}(b_{ij}+1)  = 3 \ ; \quad \sum\limits_{i=1}^n\sum\limits_{j=0}^2 N_{ij}b_{ij}=12.  
\end{eqnarray}
The above set of equations imply that
\begin{eqnarray*}
&&N_{i2} \neq 0 \implies N_{i2}=1\hspace{3pt};\hspace{3pt} N_{i0}=N_{i1}=0\\
&& N_{i2} = 0 \implies N_{i1} = N_{i0} = 1.
\end{eqnarray*}
Hence, the feasible solutions are $\{\{0,0\},\{1,0\},\{2,6\}\}$ and $\{\{0,12\},\{1,12\},\{2,0\}\}$. In order to find the correct solution we compute the flat points approximately using the discretised surface obtained from the Surface Evolver.

\subsection{Numerical computation of the flat points}\label{subsec:fps_qtz}
The discretised surface consists of vertices, edges and planar faces. The gaussian curvature at a vertex can be estimated by calculating the excess or deficit angle at that point. For a flat surface, the interior angle surrounding every vertex is  exactly $2\pi.$
When the total interior angle at vertex is more than $2\pi$, that vertex is said to be hyperbolic or a point of negative gaussian curvature. Likewise, vertices containing less than $2\pi$ would be considered spherical or be points of positive gaussian curvature. However, gaussian curvature has dimension of inverse length squared. Therefore, we estimate the gaussian curvature of the surface at a point as the excess or deficit angle divided by the total area of all the triangles concurrent at a common vertex. The accuracy of this curvature estimation increases as the resolution of the discretisation is increased. The curvature formula that we apply on the discretised surface bears a close resemblance to the continuum expression (for the gaussian curvature at a point), 
\begin{equation*}
K = \lim_{r\to 0}\frac{3}{\pi}\left(\frac{2\pi r - L_r}{r^3}\right)
\end{equation*}
where $L_r$ is the arc length of the image of a circle of radius $r$ centered at the origin in the tangent space at a point on the surface under the \emph{exponential map} \cite{doCarmo_1992}. We are interested in flat points, which have zero gaussian curvature. Thus, we look for vertices where $\vert K\vert<\varepsilon$. However, the finer the discretisation, the more vertices satisfy this condition. In order to find the actual zeros of gaussian curvature we turn to the topology of the Gauss map. Such flat points correspond to branch points in the branched cover of $S^2$. Therefore, we can use winding number around the flat point to pin down which vertices are actual flat points and which have a small, but non-zero magnitude of gaussian curvature. We compute the discrete winding number of the Gauss map around in a small neighbourhood around each of our test vertices. If the winding number exceeds one then inside this neighbourhood lies a flat point. We can close in on the flat point further by decreasing the radius of the neighbourhood. We repeat this process until we have found all the flat points to a reasonable accuracy.

This method not only allows us to find the flat points but also to classify them according to Fogden \& Hyde's regularity condition \cite{Fogden_i_1992}. We find that the winding number is precisely two for every flat point implying all the branch points are of order one. Therefore the correct branch point structure is $\{\{0,12\},\{1,12\},\{2,0\}\}$. Although this procedure does not yield accurate locations of the flat points of the QTZ-QZD surface, it reveals the correct branch point structure of the Weierstrass function without fail. There are twelve flat points mapped to twelve branch points of order one. Since there is an unbranched sheet going over every branch point the QTZ-QZD surface is classified as an \emph{irregular} surface. The branch points are given by 
\begin{equation*}
\left\{(1-\delta)e^{i\left(\frac{(2m+1)\pi}{6}+\eta\right)}\right\}_{m=0}^5\hspace{5pt};\hspace{5pt}\left\{\frac{1}{(1-\delta)}e^{i\left(\frac{(2m+1)\pi}{6}-\eta\right)}\right\}_{m=0}^5
\end{equation*}
where $-\pi/6<\eta<\pi/6$ and $0<\delta<1$. We know from the surfaces generated using Surface Evolver that the parameters $\eta$ and $\delta$ are not independent and that they are functions of the ratio of crystallographic axes or the  \emph{dimensionless chiral pitch} $\rho = c/a,$ where $c$ is the vertical period and $a$ is the horizontal period of the space group lattice. We compute the values of $(\delta,\eta)$ by numerically solving the period problem to get a family of embedded triply-periodic minimal surfaces.

\section{The Weierstrass function for the QTZ-QZD family}\label{sec:wf}
Substituting $s=3$ in \ref{eqn:qs} we get
\begin{equation}\label{eqn:qtz_qs}
\sum\limits_{i=1}^{12} q_{3,i}  = 12.
\end{equation}
And from \ref{eqn:asmp_anls} it follows that $\sum_{i=1}^{12} q_{2,i} + \textrm{deg}(P_2) \leq 8; \ \sum_{i=1}^{12} q_{1,i} + \textrm{deg}(P_1) \leq 4$.
Since $q_{3,i}$ is a non-zero positive integer and multiplicity of $\xi_i$ is one, \ref{eqn:qtz_qs} implies that $q_{3,i}=1$ for all $i\in\{1,2,...,12\}$.
Finally from \ref{eqn:bp_anls} we have $q_{0,i} = q_{1,i}  = 0 \ ; \hspace{3pt} q_{2,i}  > 0.$
Summarising all the results we have obtained so far:
\begin{equation}\label{eqn:P_pol}
  \begin{split}
 \alpha_2 = &0\\
 a_0(\xi) &= \alpha_0
 \end{split}
 \quad\quad
 \begin{split}
 q_{1,i} = 0\\
 \textrm{deg}(P_1)& = 4
 \end{split}
 \quad\quad
 \begin{split}
 a_3(\xi) = &\alpha_3\prod_{i=1}^{12}(\xi-\xi_i)
 \end{split}
\end{equation}

\subsection{The discriminant zero points}
\begin{definition}
$\xi\in\mathbb{C}$ is a discriminant zero point if and only if the discriminant of the polynomial equation \ref{eqn:pole} evaluated at $\xi$ vanishes. This implies that the roots of the polynomial evaluated at $\xi$ are not distinct.  
\end{definition}

The discriminant zero equation for a cubic polynomial, $\sum\limits_{m=0}^3\phi_{3-m}z^m=0$ is given by
 \begin{equation}\label{eqn:des}
 \phi_0(4\phi_0\phi_2^3-\phi_1^2\phi_2^2+4\phi_1^3\phi_3+27\phi_0^2\phi_3^2-18\phi_0\phi_1\phi_2\phi_3) = 0.
 \end{equation}
The zeros of  $\phi_0(\xi)$ are the branch points. Notice that the second factor is nontrivial or not identically zero if and only if the triply-periodic minimal surface is irregular. We know from section \ref{sec:fps} that the QTZ-QZD surface is an irregular surface. We use the existence of discriminant zero points other than the branch points to compute $P_1(\xi)$. To that end, we need to introduce the concept of the \emph{orbit} of a point on the complex plane under the group action.

Under the group action of the symmetry transformations an arbitrary point in the complex plane gets mapped to maximum of twelve distinct points including the original point itself. We call this set of points the {orbit} of the point of interest. It again follows from symmetry arguments that all the points within an orbit share the same orbit. 
For the group $P{6_222}$, rotations about one of the axes in the $xy$-plane produce orbits with six points, otherwise the orbits have two points, corresponding to the rotation symmetry about one of the vertical screw axes. According to \ref{eqn:P_pol} $\textrm{deg}(P_1(\xi))$ is equal to four, so the number of zeros of $P_1(\xi)$ is four. Therefore, the zeros of $P_1(\xi)$ must be $\xi_{1,1} = 0,\hspace{3pt} \xi_{1,2} = 0, \hspace{3pt}\xi_{1,3} = \infty, \hspace{3pt}\xi_{1,4} = \infty.$ 
implying $P_1(\xi) = \xi^2$. We verify below that this choice is consistent with the symmetry requirements $\xi^4a_1\left(\frac{e^{-2i\phi}}{\xi}\right)  = \xi^4\alpha_1\frac{e^{-4i\phi}}{\xi^2} = \alpha_1e^{-4i\phi}\xi^2 = e^{-4i\phi}a_1(\xi).$
The second factor of the discriminant zero equation \ref{eqn:des} can be simplified to
\begin{equation}\label{eqn:des_qtz}
4\alpha_1^3\xi^6+27\alpha_0^2\alpha_3\prod_{i=1}^{12}(\xi-\xi_i)=0.
\end{equation}
This is a degree six equation and the six roots constitute the orbit of a fixed point. Therefore, the discriminant zero points must be either $\{e^{(2m+1)i\pi/6}\}_{m=0}^5$ or $\{e^{i(m\pi/3}\}_{m=0}^5$. We observe that under the rotation fixing $(x(1),y(1),z(1))_1$, we get $(x(1),y(1),z(1))_2\neq(x(1),y(1),z(1))_3$.
Likewise, the rotation fixing $(x(1),y(1),z(1))_2$ yields $(x(1),y(1),z(1))_1\neq(x(1),y(1),z(1))_3.$
Combining these two we get $R(1)_1 \neq R(1)_2 \neq R(1)_3$. However, for the rotation fixing $(x(i),y(i),z(i))_1$, we find $(x(i),y(i),z(i))_2 = (x(i),y(i),z(i))_3$. Together, these imply $R(i)_2 = R(i)_3.$
The subscripts above are used to distinguish between sheets in the branched cover. Since the Weierstrass-Enneper representation is a bijective mapping between the branched cover and $\mathcal{M}_{\Lambda}$, the discriminant zero points are $\{e^{(2m+1)i\pi/6}\}_{m=0}^5$. Substituting $\xi=i$ in \ref{eqn:des_qtz} gives
\begin{equation}
-4\alpha_1^3+27\alpha_0^2\alpha_3\prod_{j=1}^{12}(i-\xi_j) = 0\implies\beta_1^3= \beta_0^2\alpha_3\prod_{j=1}^{12}(i-\xi_j).
\end{equation}
where $\alpha_0 = \pm2\beta_0$ and $\alpha_1 = 3\beta_1$.
\footnote{Note, a triply-periodic minimal surface from Schoen's H'-T family of surfaces \cite{Schoen_1970,Fogden_iii_1993,Fogden_1996} with space group symmetry $P\frac{6}{m}\frac{2}{m}\frac{2}{m}$ ($m$ denotes mirror planes) has the same branch point structure (including the values of the branch points) as that of the QTZ-QZD surface for $\eta = 0$ (where $\eta$ is the parameter defined in section \ref{sec:fps}\ref{subsec:fps_qtz}). The mirror symmetries of H'-T surface imply that the constants $\alpha_i$ are all real numbers and the set of nontrivial discriminant zero points is $\{e^{i(m\pi/3)}\}_{m=0}^5$. 
}

\subsection{The solution function}
The final polynomial equation for the Weierstrass function of the QTZ-QZD surface is given by
\begin{equation}\label{eqn:wf}
R^3+3\zeta^2Q(\xi)\xi^2R - 2\zeta^3Q(\xi) = 0
\end{equation}
where $\zeta = \beta_0/\beta_1$\footnote{Without loss of generality we have chosen $\alpha_0  = -2\beta_0$.} and $Q(\xi)=\prod_{i=1}^{12}(i-\xi_i)/(\xi-\xi_i)$. For each $\xi\in\mathbb{C}\cup\{\infty\}$ equation \ref{eqn:wf} is a degree three polynomial in $R$. The solution is given by
\begin{equation*}
R(\xi) = \zeta Q(\xi)^{1/3}\left[z_1\left(1+\sqrt{1+\xi^6Q(\xi)}\right)^{1/3}-z_2\left(-1+\sqrt{1+\xi^6Q(\xi)}\right)^{1/3}\right].
\end{equation*}
where $z_1$, $z_2$ are a pair of cube roots of unity satisfying $z_1z_2 = 1.$
The solution can be further simplified to
\begin{equation}
R(\xi) = z_1R_1(\xi) - z_2R_2(\xi)
\end{equation}
where
\begin{align*}
R_1(\xi) = \frac{\zeta}{\sqrt{g(\xi)}}\left(g(i)(1+\xi^6+\sqrt{g(\xi)})\right)^{1/3}; \quad R_2(\xi) = \frac{\zeta}{\sqrt{g(\xi)}}\left(g(i)(1+\xi^6-\sqrt{g(\xi)})\right)^{1/3}.
\end{align*}
and $g(\xi)=\prod_{i=1}^{12}(\xi-\xi_i)$ with $g(1)+g(i)=4$. 
Therefore, the Weierstrass function is of the form
\begin{align}\label{eqn:soln}
&R_\textrm{u}(\xi) = R_1(\xi)-R_2(\xi)\nonumber\\
&R_{\textrm{b}_1}(\xi) = e^{i(2\pi/3)}R_1(\xi)- e^{i(4\pi/3)}R_2(\xi)\nonumber\\
&R_{\textrm{b}_2}(\xi) = e^{i(4\pi/3)}R_1(\xi)- e^{i(2\pi/3)}R_2(\xi).
\end{align}
corresponding to three branches in the branched cover or the Riemann surface of interest. The labels \emph{branched} (denoted by "R$_\textrm{b}$") and \emph{unbranched} ("R$_\textrm{u}$") concern only the branch points, for e.g. $R_\textrm{u}(\xi)$ evaluated at $\xi=\xi_i$ is finite for all $i\in\{1,2,...,12\}$. There is no unique way to demarcate a boundary that separates the unbranched and branched sheets. A set of branch cuts chosen arbitrarily up to branch points 
 can be thought of as virtual boundaries between different sheets. The finite valued solution of the polynomial equation \ref{eqn:wf} at the branch points $\xi_i$ gives
\begin{equation}\label{eqn:reg_soln}
R(\xi_i) = \frac{-2\zeta}{3\xi_i^2}.
\end{equation}
It can be shown that $R_\textrm{u}$ evaluates to \ref{eqn:reg_soln} and $R_{\textrm{b}_1}$, $R_{\textrm{b}_2}$ are singular at every branch point.

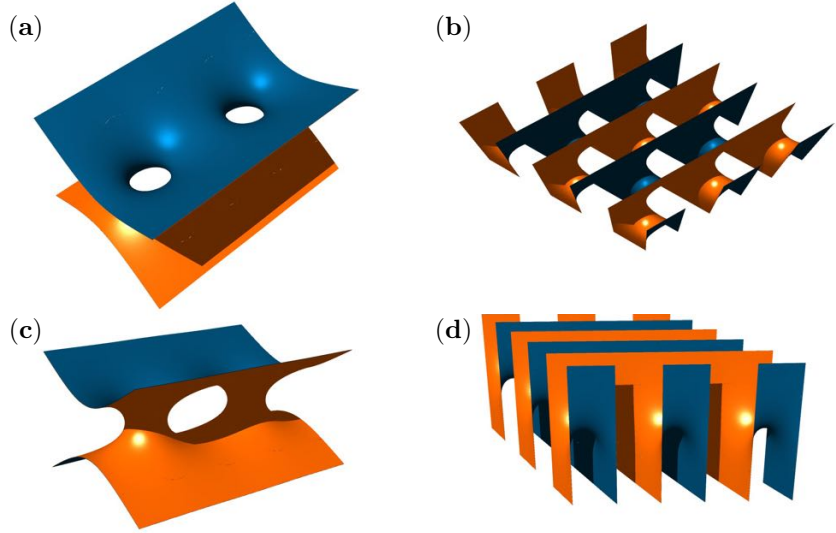
\begin{figure}[h!]
\labellist
\small\hair 2pt
\pinlabel $\bf{(a)}$ at 10 380
\pinlabel $\bf{(b)}$ at 320 380
\pinlabel $\bf{(c)}$ at 10 160
\pinlabel $\bf{(d)}$ at 320 160
\endlabellist
\includegraphics[scale=0.53]{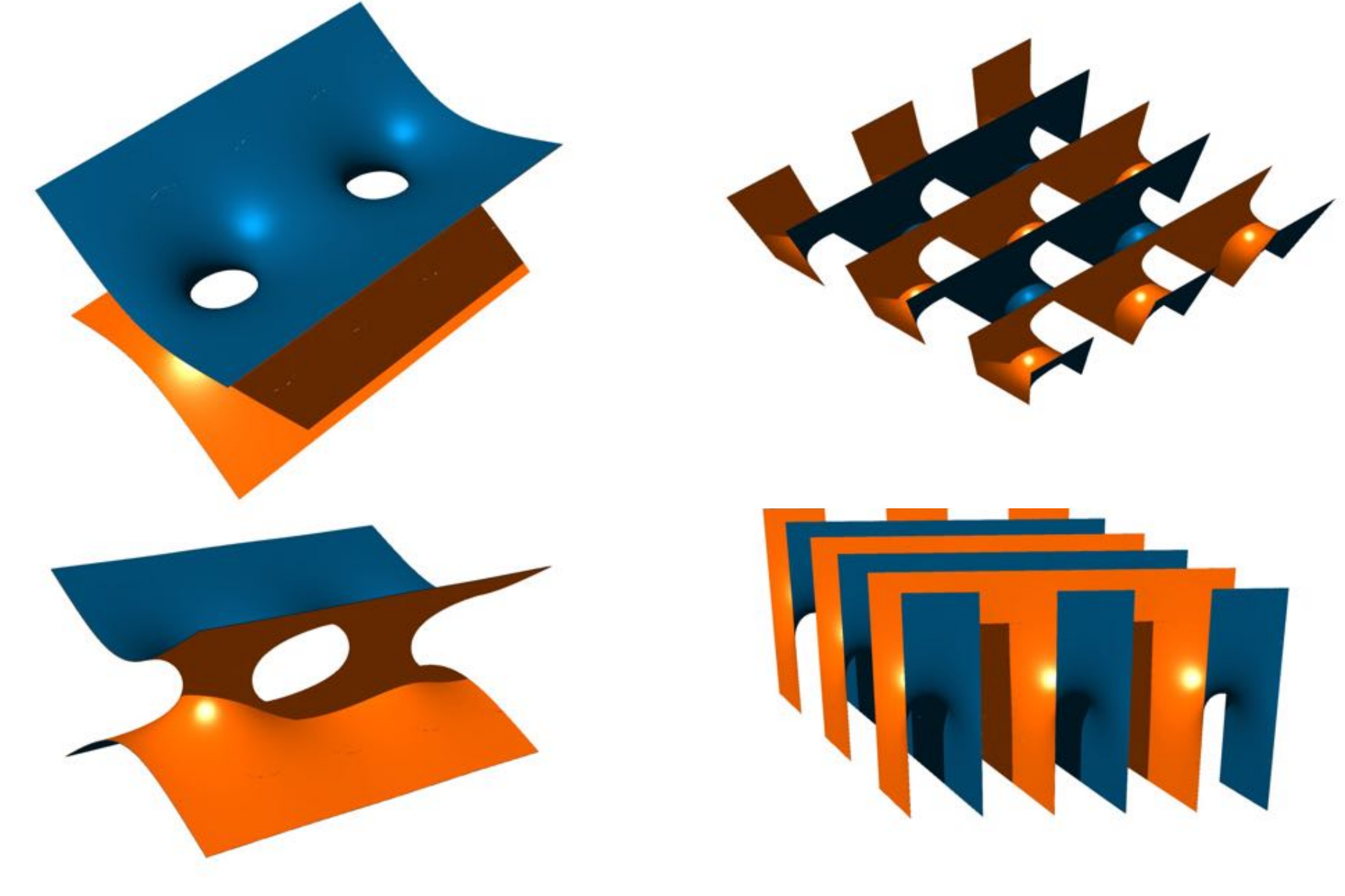}
\centering
\caption{The limiting members of the QTZ-QZD family of surfaces. (a) A top view and (b) a side view of singly periodic Scherk surface. (c) A top view and (d) a side view of doubly periodic Scherk surface.}
\label{fig:test5}
\end{figure}

The limiting members of the QTZ-QZD family of surfaces are the sheared singly periodic and the sheared doubly periodic Scherk surfaces \cite{Scherk_1835}. To see this substitute $\{e^{i(2m+1)(\pi/3)}\}_{m=0}^5$ for the branch points (with multiplicity two) in $R_\textrm{u}(\xi)$ from \ref{eqn:soln}.
\begin{align}
R_\textrm{l}^{+}(\xi) &= 4^{1/3}\left[\frac{((\xi^6+1)+(\xi^6-1))^{1/3}-((\xi^6+1)-(\xi^6-1))^{1/3}}{(\xi^6-1)}\right]\nonumber\\
& = \frac{2}{\xi^4+\xi^2+1} = \frac{2}{(\xi^2-e^{2i\pi/3})(\xi^2-e^{-2i\pi/3})}\nonumber\\
R_\textrm{l}^{-}(\xi) & = \frac{2i}{\xi^4+\xi^2+1}
\end{align}
The functions $R_\textrm{l}^{+}(\xi)$ and $R_\textrm{l}^{-}(\xi)$ are the Weierstrass functions for the doubly periodic sheared Scherk surface and the singly periodic sheared Scherk surface respectively\footnote{In this case the angle between the two ends of the singly periodic sheared Scherk surface extending to $\infty$ is $2\pi/3$.}. If we substitute $\{e^{i(2m+1)(\pi/3)}\}_{m=0}^5$ for the branch points in $R_{\textrm{b}_1}(\xi)$ and $R_{\textrm{b}_2}(\xi)$ we get functions that give rise to surfaces rotated by $-2\pi/3$ and $2\pi/3$ (with respect to the surfaces obtained from $R_\textrm{l}^{+}$ and $R_\textrm{l}^{-}$) respectively.  The minimal surfaces generated using $R_\textrm{l}^{+}(\xi)$ and $R_\textrm{l}^{-}(\xi)$ are shown in Figure 2 of the ESM. Note that the limiting members of the QTZ-QZD family of surfaces are \emph{conjugate} or \emph{adjoint} to each other. If at the very beginning of this section, we choses  $+2\zeta^3$ in \ref{eqn:wf} then the solution function is $-R$. The Weierstrass-Enneper representation with $-R$ parametrises an inversion of the QTZ-QZD surface about the origin. 
It can be shown that inverting a chiral surface about any point leads to its chiral counterpart. Therefore $-R$ is the Weierstrass function for the QTZ-QZD surface with left handed screw axes symmetries.

\begin{figure}[h!]
\labellist
\small\hair 2pt
\pinlabel $N$ at 370 150
\pinlabel $\textrm{equator}$ at 650 20
\pinlabel $(1,0,0)$ at 420 50
\endlabellist
\includegraphics[scale=0.55]{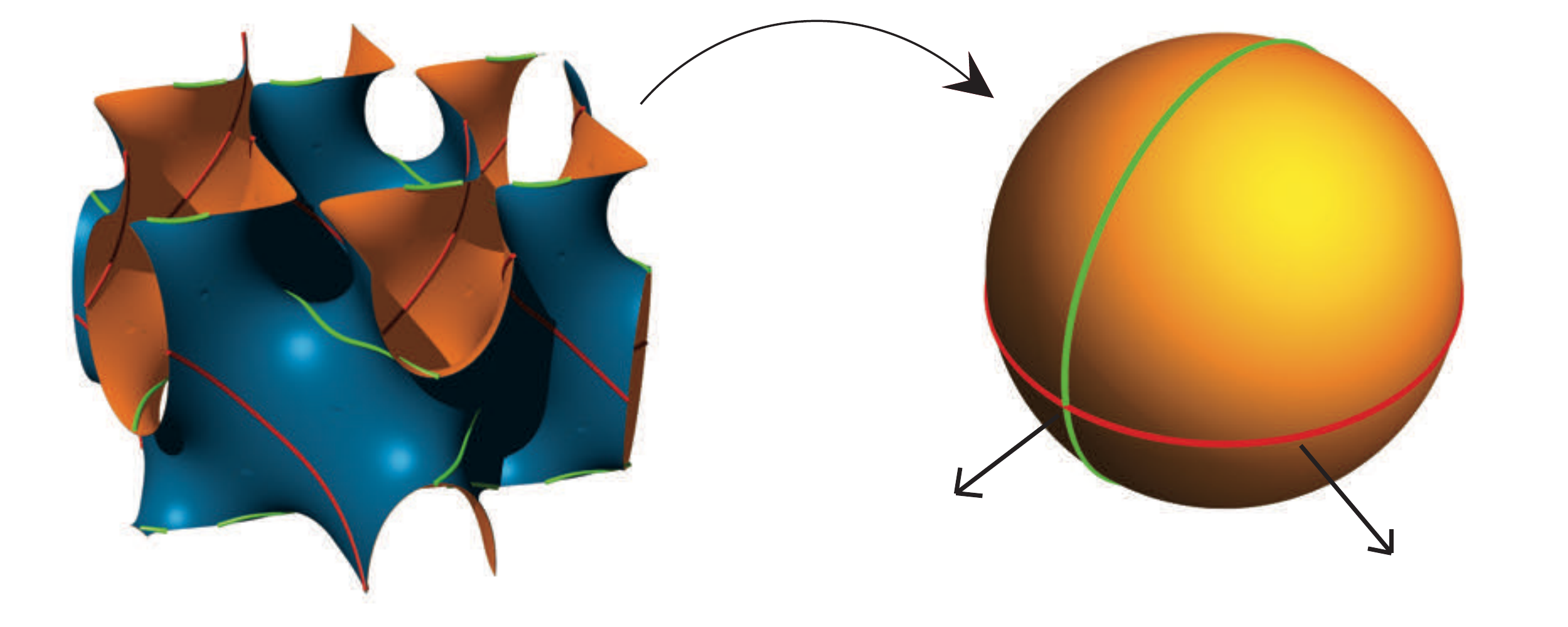}
\centering
\caption{The equator and the geodesic passing through the north pole and the point (1,0,0) on the unit sphere are mapped to the red and green curves on the QTZ-QZD surface. These two closed curves are elements of the fundamental group of $\mathcal{M}/\Lambda$ (non-contractible loops on the minimal surface). The Weierstrass-Enneper integral along the red curve gives the vertical period $c$ and the integral along green curve evaluates to the horizontal period $a$.}
\label{fig:test5}
\end{figure}

\subsection{The period problem for the QTZ-QZD family}
The Weierstrass-Enneper representation guarantees that the surface generated from \ref{eqn:wer} is minimal. Therefore, a function that we obtain by solving \ref{eqn:wf} gives us a \emph{minimal immersion} in $\mathbb{R}^3$, at every point on the surface the tangent space is spanned by two linearly independent tangent vectors and the mean curvature is zero identically. But in general this may not be triply-periodic. The period problem is nothing but solving for a two-dimensional minimal immersion into $E^3$ that is triply-periodic and hence consistent with a given set of space group symmetries and branch point structure. The domain space for the Weierstrass function is a genus four Riemann surface and hence \emph{multiply connected}. This implies the integral between any two points is multivalued if the part of the complex plane enclosed by a loop formed by two distinct paths between these two points is not simply connected or encloses branch point. In our case --  a branched cover of three sheets, the difference between multiple values of integrals (for the same set of endpoints in the complex plane) is equal to an integer multiple of one-third of the lattice parameter $c$, and this criterion solves the period problem 
by uniquely associating a numerical value for the dimensionless chiral pitch $\rho$ to each surface in the family.

The integration domain in \ref{eqn:wer} for the QTZ-QZD family of surfaces is a Riemann surface homeomorphic to a four-handled torus, a closed genus for surface. The \emph{first homology group} of the translational unit of the QTZ-QZD surface (shown in Figure \ref{fig:test3}) is $\mathbb{Z}_8$. Hence, there are eight independent non-contractible closed curves embedded in the surface referred to as the \emph{generators} of this group. The Weierstrass-Enneper integral along each of these loops gives the period associated to a different symmetry element of the space group, $P6_222$. Two loops --  one spanning the vertical period and another spanning the horizontal period -- are illustrated in Figure \ref{fig:test5}. In case of a rotational symmetry we get a constraint equation by imposing the period to be zero. If the period does not vanish in these cases then the rotational symmetry leads to self-intersections. The first two equations in \ref{eqn:pp} below ensure that the integrals along the two closed paths corresponding to the imaginary and real axis on the complex plane vanish. On the surface, these curves trace paths that loop around the handles parallel to the $x$-axis and the $y$-axis respectively, in fact the intersection-free handles are due to the vanishing periods associated to the two-fold rotational symmetries:       
\begin{equation}\label{eqn:pp}
\begin{split}
x(0)  = x(i)\hspace{5pt};
\end{split}
\quad\quad
\begin{split}
y(0)  = y(1)\hspace{5pt};
\end{split}
\quad\quad
\begin{split}
z(e^{-i\pi/3})  = z(0).
\end{split}
\end{equation}
The last equation follows from a symmetry generated by a combination of the vertical screw axis rotation and a two-fold rotation about a horizontal axis. The set of integral constraint equations in \ref{eqn:pp} constitute a sufficient condition for the parametrised surface to be a solution to the period problem. The lattice parameters of a solution to \ref{eqn:pp} are given by
\begin{align}
& z(e^{-i\pi/3})-z(1) = z(0)-z(1)=\frac{c}{6} = \frac{z(i)-z(0)}{2}\\
& \sqrt{3}(|x(e^{-i\pi/3})+(y(e^{-i\pi/3})/\sqrt{3})|+x(0)) =\frac{1}{2}a,
\end{align}
where the periods $c$ and $a$ are the vertical and horizontal translations, respectively, of the space group $P6_222$, shown in Figure \ref{fig:test7}(a). 

\begin{figure}[ht!]
\hspace{-15pt}
\labellist
\small\hair 2pt
\pinlabel $\bf{(a)}$ at 20 200
\pinlabel $\bf{(b)}$ at 430 200
\pinlabel $\theta$ at 47 175
\pinlabel $\rho$ at 360 25
\pinlabel $\textrm{Im}(\omega)$ at 580 200
\pinlabel $\textrm{Re}(\omega)$ at 760 10
\pinlabel $\theta=\pi$ at 350 160
\pinlabel $\theta=\pi/2$ at 350 55
\endlabellist
\centering
\includegraphics[scale=.465]{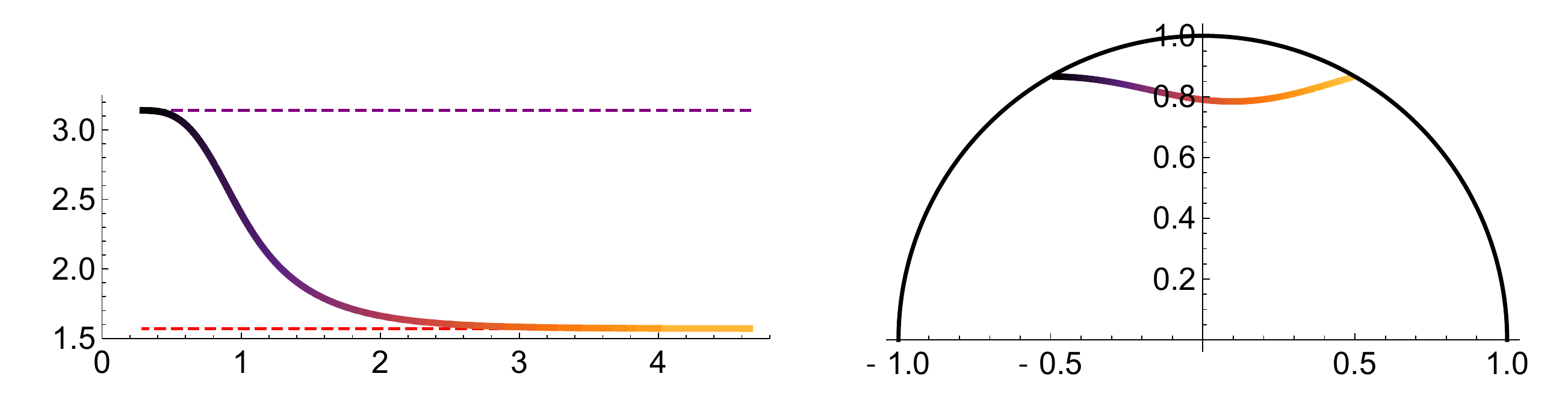}
\caption{Solution to the period problem of the QTZ-QZD family of TPMS. (a) The Bonnet angle $\theta(\rho)$ determines the parameter $\zeta$ in \ref{eqn:wf} there by completely determining the Weierstrass function for a set of branch points. (b) The branch point corresponding to $m=1$ (as per the notation used in section 4(b)). The chiral pitch $\rho$ increases monotonically with \textrm{Re}($\omega_1$).}
\label{fig:test6}
\end{figure}
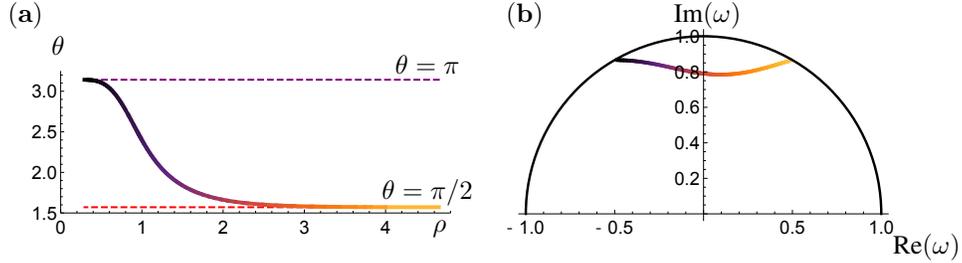

We numerically solve the period problem by varying the parameters $\delta$ and $\zeta$, whilst fixing $\eta$. Since the gaussian curvature scales identically with the magnitude of $\zeta$ without altering the global geometry, it can be reduced to $e^{i\theta}$. Hence, there are three equations involving three unknowns: $\theta$, $\delta$ and $\eta$. The equations in \ref{eqn:pp} uniquely determine a triply-periodic minimal surface up to scale and translation. Therefore, any value of ($\theta$, $\delta$ and $\eta$) satisfying \ref{eqn:pp} solves the period problem. The triply-periodic minimal surfaces thus generated form a continuous family in the dimensionless chiral pitch, $\rho$. The solution curves of the period problem are shown in Figure \ref{fig:test6}. The solution also yields $\theta_l^{+} = \pi$ and $\theta_l^{-} = \pi/2$ in agreement with the theoretical analysis in section 5(b).

\subsection{The minimal asymmetric patch}\label{sec:map}
The \emph{minimal asymmetric patch}, also known as the \emph{fundamental domain}, is defined as the smallest part of a triply-periodic minimal surface that can generate the entire surface under the repeated application of the space group elements. We obtain a minimal asymmetric patch by integrating the Weierstrass function over a suitably chosen simply connected subset of the branched cover. The minimal asymmetric patch is also called the \emph{minimal generating patch}, we shall refer to it as the minimal patch for brevity. To find the minimal patch we look for the surface patch enclosed within the primitive cell, $\mathcal{M}/\Gamma,$ where $\Gamma$ are the generators of the group  $P6_222$. 

The minimal patch plays a crucial role in constructing a triply-periodic minimal surface. 
In \cite{Koch_1988} Koch $\&$ Fischer list several triply-periodic minimal surfaces of non-cubic symmetry that are generated from the minimal patches bounded by linear asymptotes including \emph{HS1} a genus seven triply-periodic minimal surface with the space group symmetry $P6_222$. Note that this surface is distinct from the QTZ-QZD family of surfaces, as the QTZ-QZD family has genus four and lacks linear asymptotes.

\begin{figure}[h!]
\labellist
\small\hair 2pt
\pinlabel $\bf{(a)}$ at 10 300
\pinlabel $\bf{(b)}$ at 410 300
\endlabellist
\centering
\includegraphics[scale=.4]{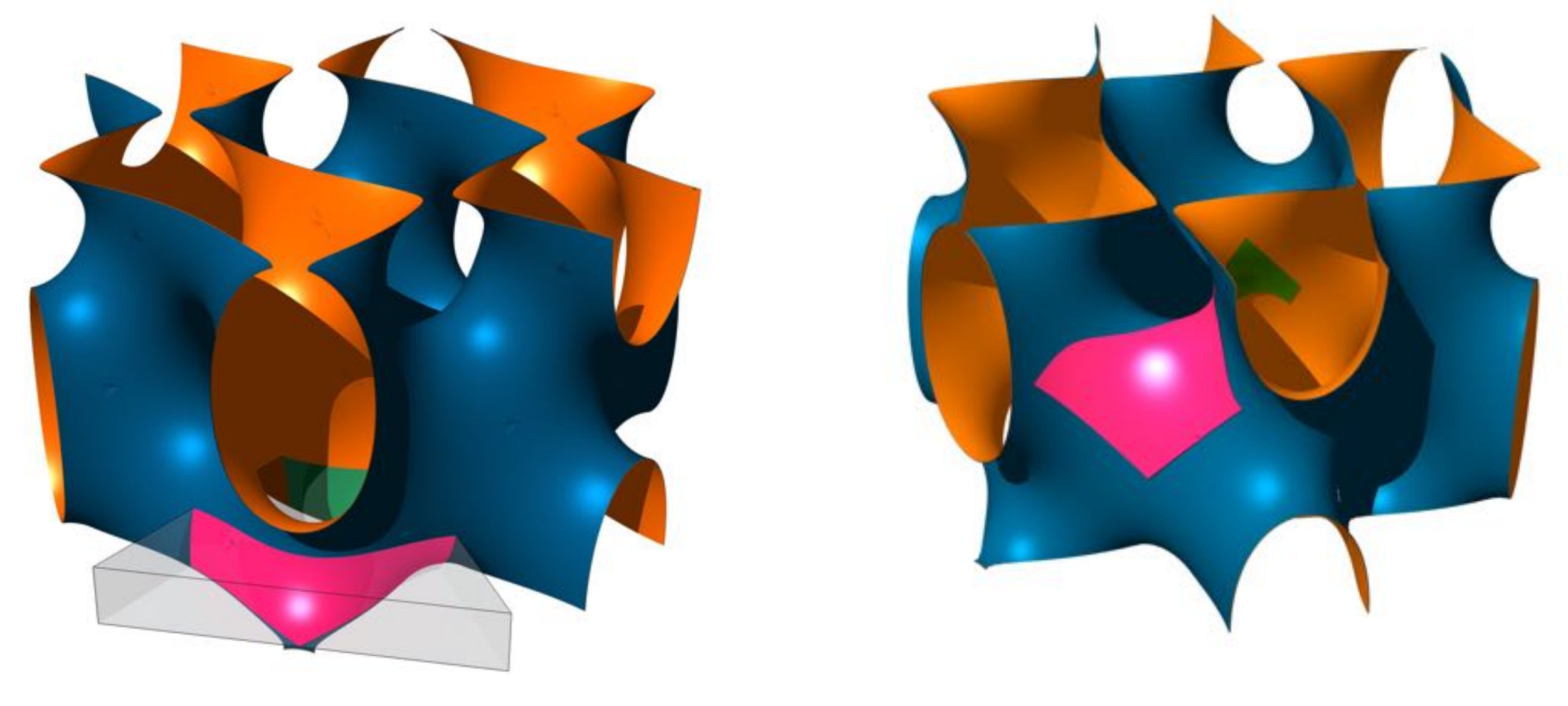}
\caption{Minimal surface patches for the QTZ-QZD surface with $\rho\approx1$ is highlighted in pink. (a) The equilateral prism shown here is a primitive cell of the translational unit cell of the $P6_222$ lattice (shown in Figure \ref{fig:test3}) and encloses a minimal patch. (b) A minimal patch bound by geodesics.}
\label{fig:test7}
\end{figure}

The surface patch within the equilateral prism shown in Figure \ref{fig:test7}(a) is a minimal patch for the QTZ-QZD surface. The prism has dimensions of $c/6$, $a$ occupying one twelfth of the volume of a translational unit cell of the lattice. Hence, the area of a minimal patch is one twelfth of the area of $\mathcal{M}_{\Lambda}$, and it contains exactly one flat point. 
Additionally, all symmetry axes incident to a minimal patch must intersect the patch only on its boundary.
In order to compute this patch via numerical integration of \ref{eqn:wer}, we need to know what part of the complex plane (domain of integration) generates it. For a discretised surface patch with arbitrary boundary curves, it would be impossible to find its exact preimage, so we look for a minimal patch with its boundaries mapped to parts of great circles on the unit sphere that pass through the images of fixed points of symmetries.

Since the image of the Gauss map of the surface forms a three sheeted branched cover of the unit sphere, a minimal patch is mapped to a region on the sphere covering an area of $\pi$. Therefore, we compute this new minimal patch by choosing an appropriate portion on the complex plane that covers an area of $\pi$ on the unit sphere under the inverse stereographic projection. A unit-radius semi-circular region centered at the origin satisfies all the essential properties: the boundary components are mapped to geodesics that pass through the fixed points of symmetries on the surface and two branched sheets are spanned (covering one flat point on the surface). A minimal patch thus obtained is shown in Figure \ref{fig:test7}(b) and the corresponding domain of integration is shown in Figure 1 of the ESM.

\section{Conclusion}
The QTZ-QZD family of surfaces is a special set of triply-periodic minimal surfaces. Like the gyroid and the lidinoid, they lack mirror symmetries and in-surface two-fold rotation axes. Yet, the absence of planar reflection symmetry implies that they are chiral as well. Since \emph{a priori} we do not know anything about the boundaries of a minimal patch, it is not possible to prove the existence using the bottom-up approach starting from its minimal patch. Instead, we have  describe a novel algorithm that proves the existence of the QTZ-QZD family of surfaces by explicit construction. Our method involves using a numeric approximation to discern the complex character of the Gauss map coupled with numeric solutions to the period problem to show that this family of surfaces has a single continuous degree of freedom related to the ratio of the lattice periods referred to as the dimensionless chiral pitch.

Particularly, we propose an algorithm that can generate triply-periodic minimal surfaces based on their skeletal graphs. We find a parametrisation (a smooth map from $\mathbb{C}$ to $\mathbb{R}^3$) based on the Weierstrass-Enneper representation. We further solve the period problem to create an embedding. The techniques used here are general and the triply-periodic minimal surfaces with in-surface symmetries can be treated as special cases.

The relevance of this work (beyond its mathematical impact) is underpinned by the fact that triply-periodic minimal surfaces already have a plethora of applications in the natural sciences and engineering disciplines \cite{Saranathan_2010,Michielsen_2008,Gerd_2011,Wilts_2017,Galusha_2008,Wilts_2012,Pendry_2004,Lu_2013,Cui_2010}: expanding the repository of known structures, while at the same time contributing to a thorough understanding of the former, can be influential in developing accurate laboratory models for experimentation.
The gyroid, for example, lends itself to self-assembly and has its own story of success in the photonic crystal and metamaterial communities \cite{Lu_2013}. Chiral metamaterials more generally provide a route to negative refraction \cite{Pendry_2004}, a key ingredient for next generation lenses and cloaking devices \cite{Cui_2010}. In this context, the new hexagonal QTZ-QZD family with a freely tunable pitch provides a new opportunity to engineer on-chip chiral photonic devices such as beamsplitters \cite{Zhang_2017,Turner_2013}, optical rotators \cite{Kim_2014} or chemical sensors \cite{Zhao_2017}.

Photonic crystals based upon the QTZ-QZD surface will further give rise to a number of deterministically induced Weyl points in the photonic band structure: these are symmetry induced band degeneracies with an associated non-trivial topological invariant that leads to unique and exotic phases of light that has been recently discovered for a chiral cubic $P2_13$ sphere packing \cite{Saba_2017aa}. The $P6_222$ symmetry of the QTZ-QZD surface is one of a handful of space groups that supports genuine deterministic Weyl points \cite{Saba_2017aa}. Further progress in this area, however, heavily relies on concrete geometries that are well understood and can be realised in practice. To this end, the analytical representation provided here can be very useful in the development of accurate triply-periodic minimal surface samples at sub-micron length scales using 3D printing and methods of self-assembly in soft matter systems such as diblock copolymers and lyotropic liquid crystals.




\section*{Appendix 1}\label{sec:apndx_i}
\subsection*{A1.1 Some basic concepts from the  differential geometry of minimal surfaces}\label{subsec:apndx_a}
In order to set up some ground work for the Weierstrass-Enneper representation we need to introduce several concepts from differential geometry of surfaces immersed in euclidean three-space. Here, we give a pedagogical overview of these concepts. This account is neither complete nor self contained but conveys the essential idea to a reader familiar with basic matrix operations and the notion of maps (that are continuous and differentiable). Consider a regular smooth surface patch, 
$\Omega$ with a parametrisation
\begin{align*}
x & :\mathcal{D}\subset\mathbb{R}^2\rightarrow\mathbb{R}^3\\
(u,v)& \mapsto(x_1,x_2,x_3)
\end{align*}
The tangent plane at $p\in\Omega$ is denoted by $T_{p}\Omega$. We choose a basis 
\begin{equation*}
\mathbf{B} = \{\partial_u x|_p,\partial_v x|_p\}
\end{equation*}
for $T_p\Omega$ so that a vector lying on the tangent plane can be expressed in this basis in the form of a column vector. 
The inner product on the vector space $\mathbb{R}^3$  i.e. the standard dot product between two vectors is denoted by 
\begin{equation*}
\langle w_1,w_2\rangle = w_1\cdot w_2.
\end{equation*} 
Now we have a naturally induced inner product in $T_p\Omega$ by the standard product in $\mathbb{R}^3$. We define the \emph{first fundamental form} in the following way:
\begin{align*}
I_p & : T_p\Omega \times T_p\Omega  \rightarrow \mathbb{R}\\
(w_1,w_2) & \mapsto w_1\cdot w_2 = w_1^t\hat{\RN{1}}_p w_2 
\end{align*}
where $\hat{\RN{1}}_p$ is the following matrix 
\begin{align*}
\hat{\RN{1}}_{p} & = \begin{bmatrix} \partial_ux\cdot\partial_ux& \partial_ux\cdot\partial_vx\\ \partial_vx\cdot\partial_ux &\partial_vx\cdot\partial_vx\end{bmatrix}\\
& = \left[(Dx)_{p}(u,v)\right]^t\left[(Dx)_{p}(u,v)\right]\\
& = g(u,v),
\end{align*}
$(Dx)_{p}(u,v)$ is the derivative matrix, $w^t$ is the transpose of $w$ and $g(u,v)$ is the \emph{metric} on $\Omega$ induced by the standard euclidean metric on $\mathbb{R}^3$. The line and area elements and other intrinsic properties of the surface are determined by the metric. The mean curvature is extrinsic, as it depends on how a two-manifold sits in a higher dimensional space (in this case $\mathbb{R}^3$). In order to define the curvatures, intrinsic and extrinsic at a point $p\in\Omega$ we need to look at the rate of change of normal vectors along the basis vectors of $T_p\Omega$. The unit normal vector at a point $p\in\Omega$ is given by
\begin{equation*}
n(p) = \frac{\partial_u x\times\partial_v x}{||\partial_u x\times\partial_v x||}.
\end{equation*}
We define the \emph{Gauss map} as the unit normal vector field on $\Omega$, 
\begin{align*}
N &: \Omega\rightarrow S^2\\
p \mapsto & n(p)\hspace{3pt}\textrm{based at}\hspace{3pt}p, 
\end{align*}
then it follows that
\begin{align*}
dN &: T_p\Omega \rightarrow T_{N(p)}S^2\\
w\mapsto & (DN)_p w\in T_{N(p)}S^2 \cong T_p\Omega  
\end{align*} 
The \emph{Shape operator}, also known as \emph{Weingarten map}, is defined as
\begin{align*}
S_{p} & : T_{p}\Omega \rightarrow T_{p}\Omega\\
 w& \mapsto -(DN)_{p}(w) = S_{p}(w)
\end{align*}  
Hence, the operator $S_{p}$ acting on $w\in T_p\Omega$ quantifies the rate of change of the unit normal vector field along $w$ at $p$.
\begin{equation*}
\hat{S}_{p} = 
\begin{bmatrix} \partial_ux\cdot S_{p}(\partial_ux)& \partial_ux\cdot S_{p}(\partial_vx)\\ \partial_vx\cdot S_{p}(\partial_ux) & \partial_vx\cdot S_{p}(\partial_vx)
\end{bmatrix}
\end{equation*}
The matrix $\hat{S}_{p}$ is real and symmetric since
\begin{equation*} 
\partial_vx\cdot S_{p}(\partial_ux)=\partial_ux\cdot S_{p}(\partial_vx).
\end{equation*} 
Therefore it is diagonalisable with real eigenvalues. The eigenvalues are called the \emph{principal curvatures}. The \emph{gaussian curvature} and the \emph{mean curvature} are then defined as
\begin{align}
K(p) & = \textrm{Det}(\hat{S}_{p})\nonumber\\
H(p) & = \textrm{Trace}(\hat{S}_{p}).
\end{align} 
For computational purposes, we recast (1.1) in a different form, to get this alternate expression we need to introduce the \emph{second fundamental form}.
\begin{align*}
\hat{\RN{2}}_{p} & : T_{p}\Omega \times T_{p}\Omega  \rightarrow \mathbb{R}\\
(w_1,w_2) & \mapsto S_{p}(w_1)\cdot w_2 = \left[\hat{S}_{p}w_1\right]^tw_2.
\end{align*}
The second fundamental form can be computed in the domain space by introducing the matrix $\hat{\RN{2}}_{p}$ given by
\begin{align*}
\hat{\RN{2}}_{p} & = \begin{bmatrix} N\cdot\partial_{uu}x& N\cdot\partial_{uv}x\\ N\cdot\partial_{vu}x & N\cdot\partial_{vv}x\end{bmatrix}\\
& = \begin{bmatrix} S_{p}(\partial_ux)\cdot\partial_ux& S_{p}(\partial_ux)\cdot\partial_vx\\ S_{p}(\partial_vx)\cdot\partial_ux &S_{p}(\partial_vx)\cdot\partial_vx\end{bmatrix}\\
& = b(u,v).
\end{align*}
Now we can express the shape operator in terms of first and second fundamental forms.
\begin{align*}
\label{eq:a1.1}
\hat{S}_{p}\hat{\RN{1}}_{p} & = \hat{\RN{2}}_{p}\nonumber\\
\hat{S}_{p} & = \hat{\RN{2}}_{p}\hat{\RN{1}}_{p}^{-1} 
\end{align*}
Taking the determinant on both sides of the above equation we can rewrite the curvatures as
\begin{eqnarray}\label{eq:so_curvatures}
K(u,v) & = \frac{\textrm{Det}(b)}{\textrm{Det}(g)}\nonumber\\
H(u,v) & = \frac{1}{2}\frac{g_{22}N\cdot\partial_{vv}x-2g_{12}N\cdot\partial_{uv}x+g_{11}N\cdot\partial_{uu}x}{\textrm{Det}(g)}.
\end{eqnarray}
A regular parametrisation, $\textrm{Det}(g)\neq0$, is said to be \emph{isothermal} or \emph{conformal} when the induced metric $g$ satisfies
\begin{equation}
g_{ij} = \lambda\delta_{ij}
\end{equation}
where $\lambda\equiv\lambda(u,v)$. A \emph{minimal surface} is a critical point of the area functional with a closed curve fixed in space as its boundary  \cite{Meeks_2012}. This is equivalent to the condition that the mean curvature vanishes identically everywhere on the surface. It follows from equation (\ref{eq:so_curvatures}) that the mean curvature vanishes if and only if
\begin{equation*}
g_{22}N\cdot\partial_{vv}x+g_{11}N\cdot\partial_{uu}x - 2g_{12}N\cdot\partial_{uv}x = 0.
\end{equation*}
From (1.3) the above equation simplifies to 
\begin{equation}\label{eq:cauchy}
\partial_{vv}x+\partial_{uu}x = 0.
\end{equation}
Hence, the isothermal coordinates of a minimal surface are harmonic functions i.e., they satisfy the Laplace equation in two variables.

\begin{lemma}[\cite{Osserman_1986}]
Let $\mathcal{M}$ be a minimal surface. Every regular point of $\mathcal{M}$ has a neighbourhood in which there exists a re-parametrisation of $\mathcal{M}$ in terms of isothermal parameters. 
\end{lemma}

\begin{lemma}[\cite{Osserman_1986}]
Let $x(u,v)\in C^2$ define a regular surface $\mathcal{M}$ in isothermal parameters. Then the condition that the coordinate functions $x(u,v)$ are harmonic is necessary and sufficient for $\mathcal{M}$ to be a minimal surface.
\end{lemma} 
$C^2$ denotes the space of all twice differentiable functions. The proof of lemma (6.2) follows from lemma (6.1) and equation (\ref{eq:cauchy}). We define a set of complex holomorphic functions as follows
\begin{equation}\label{eq:holomorphic_condition}
\phi_k(\xi) = \frac{\partial x_k}{\partial u}-i\frac{\partial x_k}{\partial v}\hspace{5pt};\hspace{5pt}\xi=u+iv.
\end{equation}
The functions $\phi_k(\xi)$  satisfy the following conditions:
\begin{eqnarray}\label{eq:phi_conditions}
\sum\limits_{k=1}^3 \phi_k(\xi)^2 & =& g_{11}-g_{22}+2ig_{12} = 0\nonumber\\
\sum\limits_{k=1}^3 |\phi_k(\xi)|^2 & =& g_{11}+g_{22} \neq 0\nonumber\\
\frac{\partial(\phi_k(\xi)+\phi_k(\xi)^*)}{\partial u} & =& -\frac{1}{i}\frac{\partial(\phi_k(\xi)-\phi_k(\xi)^*)}{\partial v}\nonumber\\
\frac{\partial(\phi_k(\xi)+\phi_k(\xi)^*)}{\partial v} & =& \frac{1}{i}\frac{\partial(\phi_k(\xi)-\phi_k(\xi)^*)}{\partial u}.
\end{eqnarray} 
The first two equations follow from the fact that the parametrisation $x_k(u,v)$ is isothermal and regular. The coordinate functions $x_k(u,v)$ are harmonic implies that the complex functions $\phi_k(\xi)$ satisfy the Cauchy-Riemann conditions for holomorphic functions.
 
\begin{lemma}[\cite{Osserman_1986}]
Let $\mathcal{D}\subset\mathbb{C}\cup\infty$, $G(\xi)$ be a meromorphic function and $F(\xi)$ be an holomorphic function in $\mathcal{D}$ such that every pole of $G(\xi)$ of degree $m$ is also a zero of $F(\xi)$ of degree at least $2m$. Then the functions
\begin{equation}
\phi_1= \frac{1}{2}F(1-G^2),\hspace{10pt}\phi_2= \frac{i}{2}F(1+G^2),\hspace{10pt}\phi_3= FG
\end{equation}
will be holomorphic in $\mathcal{D}$ and satisfy equation (1.6). Conversely every triple of functions holomorphic in $\mathcal{D}$ satisfying equation (\ref{eq:holomorphic_condition}) can be expressed in the form equation (\ref{eq:phi_conditions}), except for $\phi_1\equiv i\phi_2$ and $\phi_3\equiv0$.
\end{lemma}

The proof of this lemma follows directly from the definition of \emph{holomorphic} and \emph{meromorphic} functions. Now we state the theorem that gives a two dimensional minimal immersion in $\mathbb{R}^3$ parametrised by a set of conformal coordinates. 

\begin{theorem}[\textbf{Generalised Weierstrass-Enneper representation \cite{Meeks_1990}}]\label{thm:wer}
Any minimal surface $\mathcal{M}\subset\mathbb{R}^3$ can be parametrised up to the set of translations in $\mathbb{R}^3$ as
\begin{align*}
\mathcal{W}& :\mathcal{D}\subset\mathbb{C}\cup\{\infty\}\rightarrow\mathcal{M}\subset\mathbb{R}^3\\
\mathcal{W}(\xi)&\mapsto(x_1,x_2,x_3)\\
x_k(\xi) & = \textrm{Re}\left(\int_{\xi} \phi_k(z) dz \right)
\end{align*}
where $\phi_k$ are given by equation (\ref{eq:holomorphic_condition}) with the properties stated in equation (\ref{eq:phi_conditions}). The functions $F$ and $G$ are obtained from $\{\phi_k\}_{k=1}^3$ using lemma (6.3). Conversely the map $\mathcal{W}$ defined as above using the set of holomorphic functions $\{\phi_k\}_{k=1}^3$ gives a conformal parametrisation of a minimal  immersion in $\mathbb{R}^3$. If $\mathcal{M}$ is triply-periodic with generators of the translation $\Lambda$ and $(T^3)_{\Lambda}\cong\mathbb{R}^3/\Lambda$ is the translational unit cell defined by the generators of $\Lambda$, then $x_k(\xi)$ parametrises the quotient surface $\mathcal{M}_{\Lambda}$ immersed in $(T^3)_{\Lambda}$ given the set of {translation vectors} $P=\{\int_{\gamma}(\phi_1,\phi_2,\phi_3)|\gamma\in H_1(\mathcal{M},\mathbb{Z})\}$ are contained in $\Lambda$. 
\end{theorem}

There is a natural way to relate the meromorphic function $G$ to the Gauss map using the stereographic projection map from the unit two sphere onto the complex plane.
\begin{align}
\sigma & : S^2-\{(0,0,1)\}\rightarrow \mathbb{C}\nonumber\\
 (N_x,N_y,N_z) & \mapsto  \frac{N_x}{1-N_z}+i\frac{N_y}{1-N_z} 
\end{align}

\begin{lemma}
The meromorphic function $G$ is the Gauss map $N$ composed with the stereographic projection $\sigma$. 
\end{lemma}

From lemma (6.4) we get
\begin{align}
\sigma^{-1} & : \mathbb{C}\rightarrow S^2-\{(0,0,1)\}\nonumber\\
\sigma^{-1}(u,v) & = \left(\frac{2\textrm{Re}(G)}{|G|^2+1},\frac{2\textrm{Im}(G)}{|G|^2+1},\frac{|G|^2-1}{|G|^2+1}\right).
\end{align}
The function $G$ is called the \emph{holomorphic Gauss map}. 
The gaussian curvature of a minimal surface can be expressed in terms of the functions $F$ and $G$ as
\begin{equation*}
K(\xi)  = -\frac{1}{|F|^2}\left(\frac{2}{|G|^2+1}\right)^4\left|\frac{dG(\xi)}{d\xi}\right|^2.
\end{equation*}
Given below is a derivation of this relation.

\begin{lemma}
The Gauss map of a minimal surface is anti-conformal.
\end{lemma}
 
This is easily shown by looking at the metric induced by the Gauss map.
\begin{align*}
(DN)_{p}(w_1)\cdot(DN)_{p}(w_2)  & = S_{p}(w_1)\cdot S_{p}(w_2)\\
& = \left[\hat{S}_{p}w_1\right]^t\hat{\RN{1}}_{p}\hat{S}_{p}w_2\\
& = w_1^t\hat{S}_{p}^t\hat{\RN{1}}_{p}\hat{S}_{p}w_2\\
& = w_1^t\hat{S}_{p}\hat{S}_{p}\hat{\RN{1}}_{p}w_2.
\end{align*}
A minimal surface has zero mean curvature. Hence, the shape operator is of the form
\begin{equation*}
\hat{S}_{p} = 
\begin{bmatrix} a & b\\ b & -a
\end{bmatrix} \implies 
\hat{S}_{p}^2 =\begin{bmatrix} a^2+b^2 & 0\\0 & a^2+b^2
\end{bmatrix} 
\end{equation*}
Therefore
\begin{equation}\label{eq:shape_operator}
(DN)_{p}(w_1)\cdot(DN)_{p}(w_2) = (a^2+b^2)w_1^t\hat{\RN{1}}_{p}w_2 = -\textrm{Det}(\hat{S}_{p})w_1.w_2.
\end{equation}

We can express the gaussian curvature in terms of $F(\xi)$ and $G(\xi)$. The parametrisation given by the theorem (6.1) is conformal. Therefore equation (\ref{eq:phi_conditions}) implies
\begin{align*}
g_{11}+g_{22} & = \sum\limits_{k=1}^3|\phi_k|^2\\
2g_{11} & = |F|^2|G|^2+\frac{1}{4}\left(|F|^2(1-G^2)(1-{G^*}^2)+|F|^2(1+G^2)(1+{G^*}^2)\right)\\
g_{11} & = \frac{|F|^2}{4}(|G|^2+1)^2.
\end{align*}
Using the form of the metric on $S^2$ in spherical coordinates we get
\begin{align*}
\textrm{Det}(g_N) & = 1-\left(\frac{|G|^2-1}{|G|^2+1}\right)^2\\
& = \frac{4|G|^2}{(1+|G|^2)^2}.
\end{align*} 
From equation (\ref{eq:shape_operator}), it follows that
\begin{eqnarray}\label{eq:gaussian_curvature}
K(\xi) & =& \textrm{Det}(\hat{S}_{p})\nonumber\\
&  = &-\frac{1}{|F|^2}\left(\frac{2}{|G|^2+1}\right)^4\left|\frac{dG(\xi)}{d\xi}\right|^2
\end{eqnarray}
where $x(\xi) = p$. We use one of the forms of Weierstrass-Enneper representation [9-10] where  
\begin{equation*}
F(\xi) = 2R(\xi)\hspace{5pt};\hspace{5pt} G(\xi) = \xi.
\end{equation*}
By substituting for $F$ and $G$ from above in equation (\ref{eq:gaussian_curvature}), we get
\begin{equation}
K(\xi) = -\frac{4}{(1+|\xi|^2)^4|R(\xi)|^2}.
\end{equation}
 
\subsection*{A1.2 The transformation of the Weierstrass function under rotation}

\subsubsection*{A1.2.1 Two-fold rotational symmetries along the vertical direction}
The transformation of the holomorphic Gauss map follows from the action of the rotation matrix on the Gauss map $N$. If we let $N = (x_n,y_n,z_n)$ then
\begin{equation*}
\left[\begin{array}{c} x_n^{\prime}(\xi) \\ y_n^{\prime}(\xi) \\ z_n^{\prime}(\xi)
\end{array} \right] = \begin{bmatrix} \cos(\theta) & \sin(\theta) & 0 \\ -\sin(\theta) & \cos(\theta) & 0 \\ 0 & 0 & 1 
\end{bmatrix} \left[ \begin{array}{c} x_n(\xi) \\ y_n(\xi) \\ z_n(\xi) \end{array} 
\right] = \left[ \begin{array}{c} x_n(\xi^{\prime}) \\ y_n(\xi^{\prime}) \\ z_n(\xi^{\prime}) \end{array} 
\right]
\end{equation*}
Substituting $\theta = \pi$ and using stereographic projection of $(x_n^{\prime},y_n^{\prime},z_n^{\prime})$, we get
\begin{align*}
\xi^{\prime} & = \frac{x_n^{\prime}}{1-z_n^{\prime}}+i\frac{y_n^{\prime}}{1-z_n^{\prime}}\\
                           & = \frac{-x_n}{1-z_n}+i\frac{-y_n}{1-z_n}\\
                           & = -\xi.
\end{align*}
Using the rotation matrix given above and assuming that the base of the axis is at $(0,0,0)\in\mathcal{M}_{\Lambda}\subset\mathbb{R}^3$ we get for the $z$-coordinate
\begin{align*} 
\textrm{Re}\left(\int_0^{\xi} 2\omega R(\omega)d\omega\right) & = \textrm{Re}\left(\int_0^{-\xi} 2\omega R(\omega)d\omega\right)\\
& = \textrm{Re}\left(\int_0^{\xi} 2\omega R(-\omega)d\omega\right)
\end{align*}
A rotational symmetry about an axis normal to the surface is preserved under the Bonnet transformation. Therefore, both the real and the imaginary parts transform the same way.
\begin{equation*}
\textrm{Im}\left(\int_0^\xi 2\omega R(\omega)d\omega\right) = \textrm{Im}\left(\int_0^\xi 2\omega R(-\omega)d\omega\right).
\end{equation*}
Combining the real and imaginary parts from the above equations, we get
\begin{equation}\label{eq:R_mirror_condition}
R(\xi) = R(-\xi).
\end{equation}
From equation (\ref{eq:asymptotic_behaviour}) in the main text and the above equation (\ref{eq:R_mirror_condition}), we conclude that the sets of branch points with identical branch point structure must remain invariant. Therefore, for the case of the QTZ-QZD surface, we get
\begin{equation*}
 \{\xi_i\}_{i=1}^{12} = \{-\xi_i\}_{i=1}^{12}.
\end{equation*} 
The invariance of the polynomial equation governing $R$ and equation (\ref{eq:R_mirror_condition}) imply
\begin{equation*}
a_m(\xi) = a_m(-\xi).
\end{equation*}
The form of $a_m(\xi)$ given by equation (\ref{eqn:pol_coef}) from the main text imposes the following constraints  
\begin{align*}
 \{\xi_{m,i}\} & = \{-\xi_{m,i}\} \\
 (-1)^{\textrm{deg}(P_m)+\sum\limits_{i=1}^n q_{m,i}} & = 1.
\end{align*}
where $\{\xi_{m,i}\}$ are the zeros of $P_m(\xi)$.

\subsubsection*{A1.2.2 Two-fold rotational symmetries along the horizontal directions}
There are in total six distinct two-fold rotation symmetries of $\mathcal{M}$ up to translation such that the symmetry axes lie on the $xy$-plane. Let us denote the direction of one such axis by a unit vector $(n_x,n_y,0)$ such that $\tan((\pi/2)-\phi) = n_y/n_x$.
The transformation of the unit normal vector $(x_n(\xi),y_n(\xi),z_n(\xi))$ is given by  
\begin{equation*}
\begin{bmatrix} \cos(2\phi) & -\sin(2\phi) & 0 \\ -\sin(2\phi) & -\cos(2\phi) & 0 \\ 0 & 0 & -1 
\end{bmatrix}
\left[ \begin{array}{c} x_n(\xi) \\ y_n(\xi) \\ z_n(\xi) \end{array} 
\right] = \left[ \begin{array}{c} x_n(\xi^{\prime}) \\ y_n(\xi^{\prime})\\ z_n(\xi^{\prime})\end{array} 
\right]
\end{equation*}
Therefore
\begin{align*}
\xi^{\prime} & = \frac{x_n^{\prime}}{1-z_n^{\prime}}+i\frac{y_n^{\prime}}{1-z_n^{\prime}}\\
                           & = \frac{x_n(\xi)\cos{2\phi}-y_n(\xi)\sin{2\phi}}{1+z_n(\xi)}+i\frac{-x_n(\xi)\sin{2\phi}-y_n(\xi)\cos{2\phi}}{1+z_n(\xi)}\\
                           & = \frac{(\cos{2\phi}-i\sin{2\phi})x_n(\xi)}{1+z_n(\xi)}+i\frac{(-\cos{2\phi}-i\sin{2\phi})y_n(\xi)}{1+z_n(\xi)}\\
                           & = \frac{\exp{(-2i\phi)}}{\xi}.
\end{align*}
The last equality follows from $(1+z)/(1-z) = |\xi|^2$ which is obtained using the fact that $(x_n(\xi),y_n(\xi),z_n(\xi))$ is a unit vector and $|\xi|^2=(x_n^2+y_n^2)/(1+z_n)^2$. Let the base of the chosen axis be at $(x_s,y_s,z_s)\in\mathcal{M}_{\Lambda}\subset\mathbb{R}^3$. Then we have
\begin{equation*}
\begin{bmatrix} \cos(2\phi) & -\sin(2\phi) & 0 \\ -\sin(2\phi) & -\cos(2\phi) & 0 \\ 0 & 0 & -1 
\end{bmatrix}
\left[ \begin{array}{c} x(\xi)-x_s \\ y(\xi)-y_s \\ z(\xi)-z_s \end{array} 
\right] = \left[ \begin{array}{c} x(\xi^{\prime})-x_s \\ y(\xi^{\prime})-y_s\\ z(\xi^{\prime})-z_s \end{array} 
\right]
\end{equation*}
The transformation of the $z$-coordinate gives 
\begin{align*} 
\textrm{Re}\left(\int_{0}^{\xi}2\omega R(\omega)d\omega\right) & = 2z_s - \textrm{Re}\left(\int_{0}^{\frac{e^{-2i\phi}}{\xi}} 2\omega R(\omega)d\omega\right)\\
\textrm{Re}\left(\int_{0}^{\xi}2\omega R(\omega)d\omega\right) & = 2 \textrm{Re}\left(\int_0^{\xi_1}2\omega R(\omega)d\omega\right)\\ 
& +\textrm{Re}\left(\int_{\infty}^{\xi}\frac{2e^{-2i\phi}}{\omega}R\left(\frac{e^{-2i\phi}}{\omega}\right)\frac{e^{-2i\phi}}{\omega^2}d\omega\right)\\
\textrm{Re}\left(\int_{\xi_1}^{\xi}2\omega R(\omega)d\omega\right) & = \textrm{Re}\left(-\int_{\infty}^{\xi_1}\frac{2e^{-2i\phi}}{\omega}R\left(\frac{e^{-2i\phi}}{\omega}\right)\frac{e^{-2i\phi}}{\omega^2}d\omega\right)\\
& + \textrm{Re}\left(\int_{\infty}^{\xi}\frac{2e^{-2i\phi}}{\omega}R\left(\frac{e^{-2i\phi}}{\omega}\right)\frac{e^{-2i\phi}}{\omega^2}d\omega\right)\\
& = \textrm{Re}\left(\int_{\xi_1}^{\xi}\frac{2e^{-2i\phi}}{\omega}R\left(\frac{e^{-2i\phi}}{\omega}\right)\frac{e^{-2i\phi}}{\omega^2}d\omega\right)
\end{align*}
where $(x(\xi_1),y(\xi_1),z(\xi_1)) = (x_s,y_s,z_s)$ and $\xi_1^2 = e^{-2i\phi}$. Since a symmetry transformation with the axis of rotation along a normal to the surface holds for all the surfaces in the associate family, we have
\begin{equation*}
\textrm{Im}\left(\int_{\xi_1}^{\xi}2\omega R(\omega)d\omega\right)  =  \textrm{Im}\left(\int_{\xi_1}^{\xi}\frac{2e^{-2i\phi}}{\omega}R\left(\frac{e^{-2i\phi}}{\omega}\right)\frac{e^{-2i\phi}}{\omega^2}d\omega\right).
\end{equation*}
Combining the real and imaginary parts from the above two equations, we  get
\begin{equation*} 
\xi^4 R(\xi) = e^{-4i\phi}R(e^{-2i\phi}/\xi).
\end{equation*}
As before, the branch points get mapped to the other branch points,
\begin{equation*}
 \{\xi_i\}_{i=1}^{12} = \{e^{-2i\phi}/\xi_i\}_{i=1}^{12}.
\end{equation*}
The polynomial equation must be invariant so that the solution function obtained generates the original surface, which implies
\begin{equation*}
\sum\limits_{m=0}^{n}a_m(\xi)R(\xi)^m = \sum\limits_{m=0}^{n}a_m(\xi^{\prime})R(\xi^{\prime})^m.
\end{equation*}
Now by equating the transformed polynomial and the original one term by term, we get a transformation rule for the coefficient functions
\begin{eqnarray}\label{eq:coef_tmns}
&&a_m(\xi)R(\xi)^m  = a_m(\xi^{\prime})R(\xi^{\prime})^m\nonumber\\
&&a_m(\xi)(\xi^{-4}e^{-4i\phi}R(e^{-2i\phi}/\xi))^m  = a_m(e^{-2i\phi}/\xi)R(e^{-2i\phi}/\xi)^m\nonumber\\\ 
&& a_m(e^{-2i\phi}/\xi)  = \xi^{-4m} e^{-4i\phi m}a_m(\xi).
\end{eqnarray} 
Substituting $e^{-2i\phi}/\xi$ as the argument of $a_m$ in equation (\ref{eqn:pol_coef}) in the main text and simplifying the left hand side of equation (\ref{eq:coef_tmns}), we get
\begin{align*}
\alpha_m & \prod_{i=1}^{\textrm{deg}(P_m)}\left(\frac{e^{-2i\phi}}{\xi}-\xi_{m,i}\right)\prod_{i=1}^{12}\left(\frac{e^{-2i\phi}}{\xi}-\xi_i\right)^{q_{m,i}} \\
 & = (-1)^{\textrm{deg}(P_m)+\sum\limits_{i=1}^{12}q_{m,i}}\left(\prod_{i=1}^{\textrm{deg}(P_m)}\xi_{m,i}\right)\left(\prod_{i=1}^{12}\xi_i^{q_{m,i}}\right) \frac{a_m(\xi)}{\xi^{\textrm{deg}(P_m)+\sum\limits_{i=1}^{12}q_{m,i}}}\\
& = \textrm{rhs of equation}\hspace{3pt}(\ref{eq:coef_tmns}) = \xi^{-4m}e^{-4i\phi m}a_m(\xi).
\end{align*}
The last equality follows if and only if
\begin{align*}
 4m & = \textrm{deg}(P_m) + \sum\limits_{i=1}^{n}q_{m,i} \nonumber\\
 1 & = (-1)^{\textrm{deg}(P_m)+\sum\limits_{i=1}^{12}q_{m,i}}\exp(4i\phi m)\left(\prod_{i=1}^{n}\xi_i^{q_{m,i}}\right)\left(\prod_{j=1}^{\textrm{deg}(P_m)}\xi_{m,i}\right)\nonumber\\
 \{\xi_{m,i}\} & = \{\exp(-2i\phi)/\xi_{m,i}\}
\end{align*}
for all $m\in\{0,1,...,s\}$. 

\subsubsection*{A1.2.3 Three-fold screw axis rotational symmetries along the vertical direction}
The screw axis rotation symmetry involves rotation by $(2\pi/3)$ about an axis parallel to the $z$-axis followed by a translation of $(c/3)$ units along the same. The $3\times3$ rotation matrix is identical to the one used in the beginning of section Appendix A1.2.1 with $\theta = 2\pi/3$. For the right hand screw axis symmetry, the Gauss map transforms as   
\begin{equation*}
\xi^{\prime}=e^{i(2\pi/3)}\xi.
\end{equation*}
The transformation of the $z$-coordinate leads to
\begin{align*} 
\textrm{Re}\left(\int_{0}^{\xi}2\omega R(\omega)d\omega\right) + (c/3) & = \textrm{Re}\left(\int_{0}^{e^{i(2\pi/3)}\widetilde{\xi}}2\omega R(\omega)d\omega\right)\\
                  & =  \textrm{Re}\left(\int_{\widetilde{0}}^{0}2e^{(4i\pi/3)}\omega R\left(e^{(2i\pi/3)}\omega\right)d\omega\right)\\
                  & + \textrm{Re}\left(\int_{0}^{\widetilde{\xi}}2e^{(4i\pi/3)}\omega R\left(e^{(2i\pi/3)}\omega\right)d\omega\right)\\
                  & = \textrm{Re}\left(\int_{\widetilde{0}}^{\widetilde{\xi}}2e^{(4i\pi/3)}\omega R\left(e^{(2i\pi/3)}\omega\right)d\omega\right)\\
                  & = \textrm{Re}\left(\int_{0}^{\xi}2e^{(4i\pi/3)}\widetilde{\omega} R\left(e^{(2i\pi/3)}\widetilde{\omega}\right)d\widetilde{\omega}\right)
\end{align*}
where $\widetilde{\xi}$ and $\xi$ both denote a single point on the complex plane but belong to two different sheets in the branched cover. In the second equality, we have used
\begin{equation*}
c/3 = \textrm{Re}\left(\int_{0}^{\widetilde{0}}2e^{(4i\pi/3)}\xi R\left(e^{(2i\pi/3)}\xi\right)d\xi\right).
\end{equation*}
Therefore, the action of the screw axis symmetry yields the following transformation of the Weierstrass function 
\begin{equation*}
R(\xi) = e^{(4i\pi/3)}R(e^{(2i\pi/3)}\widetilde{\xi}).
\end{equation*}

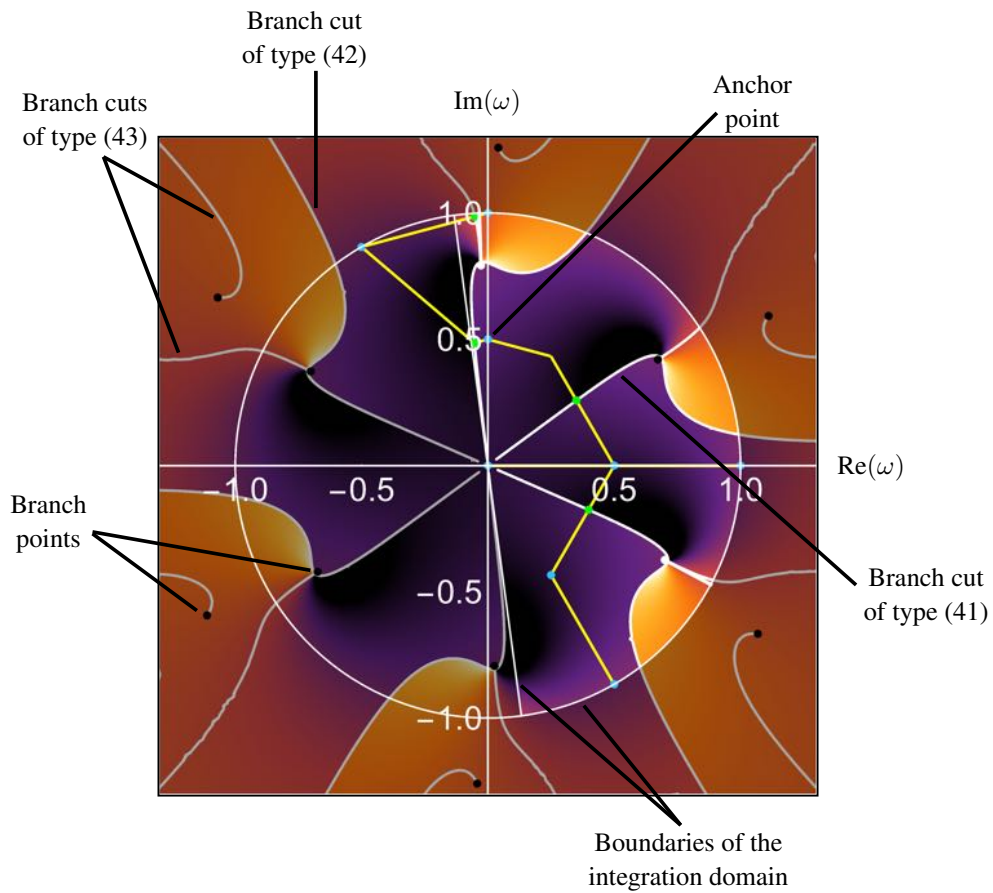
\begin{figure}
\labellist
\small\hair 2pt
\pinlabel $\textrm{Branch cut}$ at 150 495
\pinlabel $\textrm{of type (\ref{eq:branch_2})}$ at 150 475
\pinlabel $\textrm{Branch cuts}$ at 30 450
\pinlabel $\textrm{of type (\ref{eq:branch_3})}$ at 30 430
\pinlabel $\textrm{Branch}$ at 10 230
\pinlabel $\textrm{points}$ at 10 210
\pinlabel $\textrm{Anchor}$ at 395 460
\pinlabel $\textrm{point}$ at 395 440
\pinlabel $\textrm{Branch cut}$ at 490 190
\pinlabel $\textrm{of type (\ref{eq:branch_1})}$ at 490 170
\pinlabel $\textrm{Boundaries of the}$ at 360 45
\pinlabel $\textrm{integration domain}$ at 360 25
\pinlabel $\textrm{Im}(\omega)$ at 250 450
\pinlabel $\textrm{Re}(\omega)$ at 460 250
\endlabellist
\centering
\includegraphics[scale=.7]{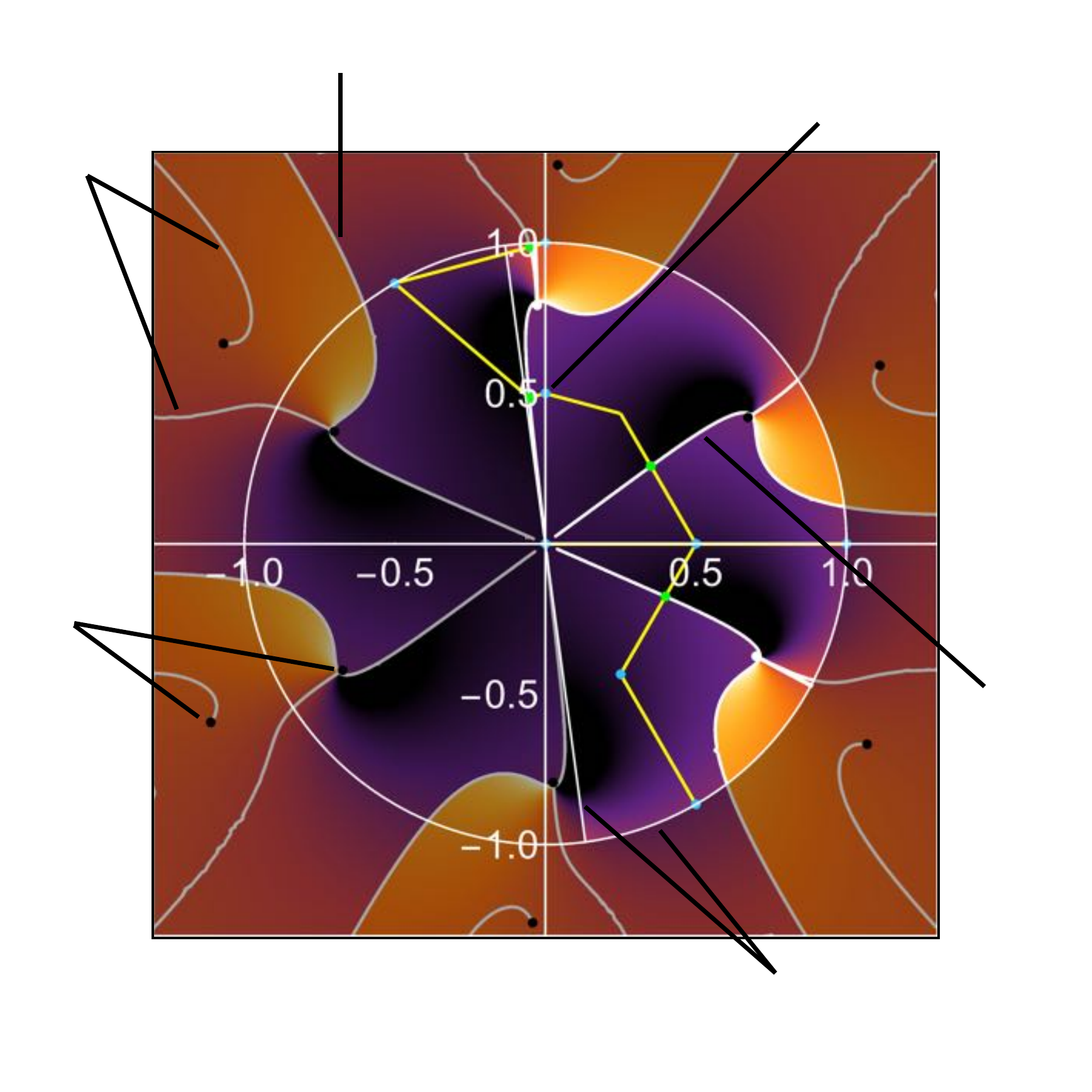}
\caption{A contour plot of the Weierstrass function of the QTZ-QZD surface for $\rho\approx 1$. The highlighted region, a semi-circle of unit radius centered at the origin is the integration domain used to generate the minimal patch shown in Figure 8(b). The piecewise linear paths (in yellow) are the integration paths used to solve the period problem. The points marked in green are the points of discontinuity of the branches of the function in (\ref{eqn:soln}) evaluated in {\tt Mathematica}.}
\label{fig:test9}
\end{figure}

\section*{Appendix 2: Numerical methods}\label{sec:apndx_ii}
\subsection*{A2.1 Branch cuts and the numerical integration}
Due to the complicated form of the Weierstrass function $R(\xi)$ in (\ref{eqn:reg_soln}), we resort to numerical integration for generating the minimal patch and solving the period problem of the QTZ-QZD family of surfaces. We use the {\tt NIntegrate} function in {\tt Mathematica}. The function is very well behaved, which implies that the integral is convergent and well defined except in the proximity of the branch points. Whilst numerically integrating a function defined on a Riemann surface, it is necessary to introduce branch cuts. Suppose we have a multivalued function, $h(\xi)^{m/n}$ such that $m$ and $n$ are co-prime then a branch cut is a\footnote{It is not unique up to the end points which must be the branch points.} continuous curve in the complex plane across which the argument of the function, $\textrm{Arg}(h(\xi)^{m/n})$ has a discontinuity of $2\pi(m/n)$. In {\tt Mathematica}, the branch cuts are placed according to a simple rule: the curve $\{\xi| \,\textrm{Im}(h(\xi))=0,\textrm{Re}(h(\xi))<0\}$ is a branch cut if $h(\xi)$ is a polynomial in $\xi$, otherwise it will apply this rule to every polynomial function constituting a compound function. Therefore for the function $R(\xi)$ the branch cuts are given by  
\begin{eqnarray}
&&\label{eq:branch_1}\{\xi\vert \, \textrm{Im}(1+\xi^6+\sqrt{g(\xi)})=0,\hspace{5pt}\textrm{Re}(1+\xi^6+\sqrt{g(\xi)})<0\}\\
&&\label{eq:branch_2}\{\xi\vert \,\textrm{Im}(1+\xi^6-\sqrt{g(\xi)})=0,\hspace{5pt}\textrm{Re}(1+\xi^6-\sqrt{g(\xi)})<0\}\\
&&\label{eq:branch_3}\{\xi\vert \,\textrm{Im}(g(\xi))=0,\hspace{5pt}\textrm{Re}(g(\xi))<0\}. 
\end{eqnarray}

\begin{itemize}
\item A branch cut where $\textrm{Im}(g(\xi))$ is zero leads to a discontinuity in $R$ such that $\sqrt{g(\xi)}$ jumps to $-\sqrt{g(\xi)}$ and $-\sqrt{g(\xi)}$ to $\sqrt{g(\xi)}$. This is because $\sqrt{g(\xi)}$ in the expression of $R(\xi)$ switches between the branches. 
\item For the branch cuts of the type equation (\ref{eq:branch_1}), $R$ jumps so that the value of $(1+\xi^6+\sqrt{g(\xi)})^{1/3}$ across the cut differs by a factor of $e^{i 2\pi/3}$ or $e^{-i 2\pi/3}$ depending on the change in $\textrm{Arg}(\xi)$ across the cut. The same is true for the branch cuts of the type equation (\ref{eq:branch_2}) where $(1+\xi^6-\sqrt{g(\xi)})^{1/3}$ is discontinuous across the cut.
\end{itemize}

Now that we know how to account for the function jumps, it would be ideal if we could compute the branch cuts analytically. But the expressions for the imaginary part of the expressions in equations (\ref{eq:branch_1}) and (\ref{eq:branch_2}) are very cumbersome to work with analytically. However, the {\tt ContourPlot} function in {\tt Mathematica} is sufficient to roughly locate the branch cuts. All the three possible branch cuts are  shown here in ESM Figure 1. When integrating close to a branch point, we require an accurate location of the branch cuts, in which case we zoom in on the region by setting the range of the {\tt ContourPlot} appropriately. Note that these plots are used only to visualise the placement of the branch cuts by an in built algorithm in {\tt Mathematica}.

While integrating numerically, we compensate on the compromise made in locating the branch cuts by controlling the path of integration minimally. Given the limits of integration, we first choose an integration path. This is done such that the roots of the polynomial equation (\ref{eqn:reg_soln}) of the main text for some $\xi\in\mathbb{C}$ on the integration path are consistent with the value of $R$ evaluated at $\xi$ on the branch traversed during the integration. We choose piecewise linear curves as paths of integration where different subpaths are joined by {anchor points} and some of the subpaths intersect the branch cuts. We can always find the subpaths that intersect the branch cuts using {\tt ContourPlot} to evaluate the real or imaginary part of $R(\xi)$. Then we can compute the point of intersection of a branch cut and a subpath by solving for the roots of equations (\ref{eq:branch_1}), (\ref{eq:branch_2}) and (\ref{eq:branch_3}) using {\tt NSolve}. The final integration path consists of these points of discontinuity with markers to ensure that during the integration the function changes the branches if encountered with a marked point according to the value of its marker.

\begin{figure}[t]
\labellist
\small\hair 2pt
\pinlabel $\bf{(a)}$ at 0 600
\pinlabel $\bf{(b)}$ at 440 600
\pinlabel $\bf{(c)}$ at 0 220
\pinlabel $\bf{(d)}$ at 440 220
\endlabellist
\centering
\includegraphics[scale=0.5]{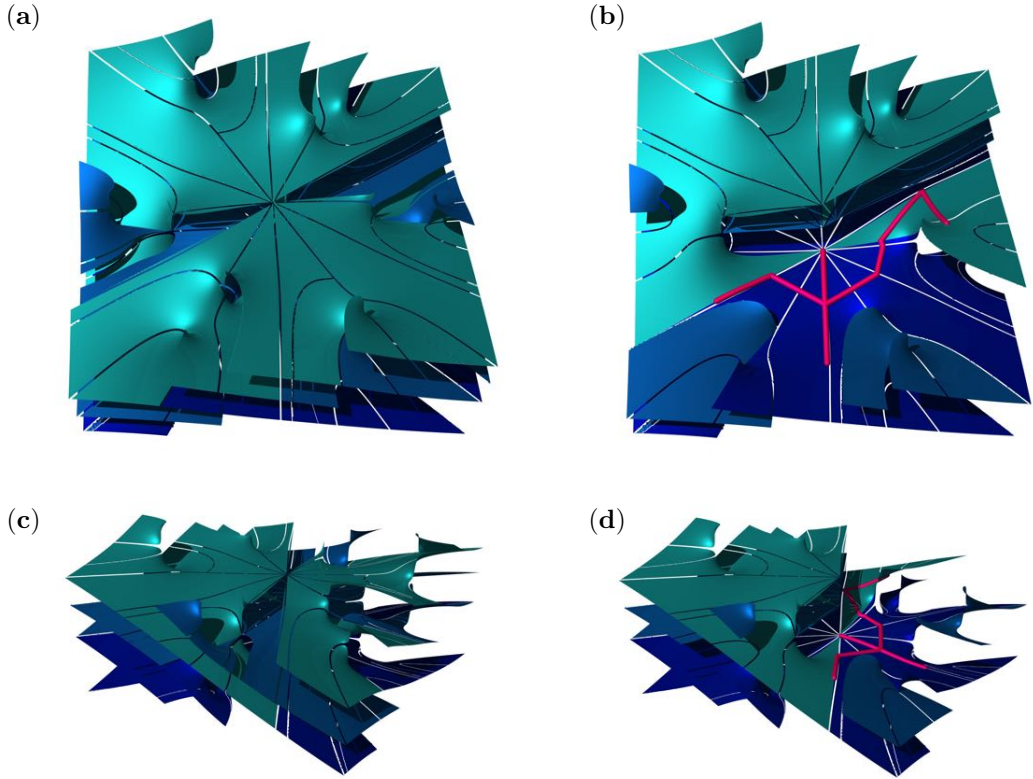}
\caption{A visualisation of the Riemann surface for the case with $\rho\approx1$. There are three sheets covering every point in $\mathbb{C}\cup\{\infty\}$ except at the branch points where two sheets are pinned. The piecewise smooth curves highlighted in pink in (b) and (d) lie on a single sheet. (a) $\&$ (b) A top view, (c) $\&$ (d) a side view. In (b) $\&$ (d) a half of the sheet on the very top is removed to show the remaining two sheets below.}
\label{fig:test4}
\end{figure}

\subsection*{A2.2 The numerical method used to construct a discretised representation of the surfaces from the QTZ-QZD family}\label{subsec:apndx_d}
Following the philosophy set out by O'Keeffe, Hyde and others, the quartz (or {qtz}) net
and the {qzd} net, shown in Figure 1(a) and (b), 
are an adequate choice for the generation of a well-defined dividing surface, as they are duals\footnote{The notion of duality refers to one in Appendix II of \cite{Schoen_1970}} of one another and of the same symmetry
(called proper duals \cite{Blatov_2007}). 
The symmetry group of both nets is the chiral hexagonal group $P6_222$ and the genus of each of the nets is $g=4$ for the portion corresponding to the primitive translational unit cell (for a net, the genus is $g=1+e-v$ if $e$ is the number of edges and $v$ is the number of vertices of the primitive cell).

As the next stage of the construction, a discrete representation (triangulation) of an interface between the {qtz} and the {qzd} nets can be generated which is neither smooth nor area-minimising. This can be the tubular surface, or a discretisation thereof, around one of the two networks; or a culled Voronoi diagram of a set of points on the two graphs (with all Voronoi faces removed that separate points on the same graph) \cite{de_Campo_2013}.
The resulting triangulated interfacial surface can be evolved towards a minimal surface, by using a conjugate gradient scheme with a function $f[S]=\int_S (H-H_0)^2 dS$ with $H_0=0$ and $H$ the point-wise mean curvature. Leaving mathematical details aside, successful minimisation of this functional to $f[S]=0$ yields the same minimal surface as does area-minimisation with a fixed (but here unknown) volume constraint \cite{KGB_1997}. For a general input surface, it is however not guaranteed that a minimum with constant $H=0$ can be reached by a topology- and/or symmetry-preserving evolution of the interface; consider as a simple example two parallel circles off-set perpendicularly by a distance $L$; for $L$ smaller than a critical value $L_c$ the unique minimal surface suspended from these circles is the catenoid, which an evolution starting from a cylindrical surface can reach. If $L>L_c$, the minimal surface consists in two flat sheets (one in each circle); however, starting from a cylindrical surface, a topology-preserving evolution under $f[S]$ will result in a pinch-off of the surface. 

For the the { qtz}/{ qzd} system, we have used Ken Brakke's conjugate gradient solver {\tt Surface Evolver} \cite{Brakke_1992} to determine a surface minimising $f[S]$, starting from a tubular surface around the {qtz} net. To a high degree of numerical accuracy, the resulting surface has $f[S]=0$, with discrete point-wise mean curvatures verified to vanish. By computer-visual inspection, all symmetries of the original symmetry group of the graphs remain valid.

\subsection*{Acknowledgements}
{
We are grateful to Michael O'Keeffe for his suggestion of the qtz-qzd as a chiral pair of networks with a tunable pitch. We are grateful to Rob Kusner and Kersten Gro\ss e-Brauckmann for early discussions of the period problem. Much of the inspiration for the work shown here stems from published concepts by Alan Schoen, Stephen Hyde, Andrew Fogden and Micheal O'Keeffe, which we here gratefully acknowledge.}



\end{document}